\documentclass[journal,  10pt, final]{IEEEtran}
\usepackage{amsmath}
\usepackage{bm}
\usepackage{amsfonts}
\usepackage{amssymb }
\usepackage{lineno}
\usepackage{graphics} 
\usepackage{mathtools}
\usepackage{cite}
\usepackage{multirow}
\usepackage{threeparttable}
\usepackage{varwidth}
\usepackage{booktabs}
\usepackage{xr}
\DeclareMathOperator*{\argmin}{arg\,min}

\usepackage{algorithmicx}
\usepackage{algorithm}
\usepackage[noend]{algpseudocode}
\algnewcommand\algorithmicinput{\textbf{Input:}}
\algnewcommand\Input{\item[\algorithmicinput]}

\algnewcommand\algorithmicoutput{\textbf{Output:}}
\algnewcommand\Output{\item[\algorithmicoutput]}

\algnewcommand{\LineComment}[1]{\Statex \(\triangledown \) #1}  

\algnewcommand{\IfThenElse}[3]{
	\State \algorithmicif\ #1\ \algorithmicthen\ #2\ \algorithmicelse\ #3 \algorithmicend}
\algnewcommand{\IfThen}[2]{\State \algorithmicif\ #1\ \algorithmicthen\ #2 \algorithmicend}

\makeatletter
\let\MYcaption\@makecaption
\makeatother

\makeatletter
\let\@makecaption\MYcaption
\makeatother

\graphicspath{{../image/, ./}}
\DeclareGraphicsExtensions{.pdf,.jpeg,.png,.tiff}

\newtheorem{definition}{Definition}
\newtheorem{assumption}{Assumption}
\newtheorem{remark}{Remark}
\newtheorem{lemma}{Lemma}
\newtheorem{proposition}{Proposition}
\newtheorem{theorem}{Theorem}
\newtheorem{problem}{Problem}
\newtheorem{example}{Example}
\newtheorem{procedure}{Procedure}

\usepackage{hyperref}
\usepackage{xcolor}
\hypersetup{
	colorlinks,
	linkcolor={red!50!black},
	citecolor={blue!50!black},
	urlcolor={blue!80!black}
}

\usepackage{pdfpages}

\newcommand{\cM}{\mathcal M}
\newcommand{\cL}{\mathcal L}
\newcommand{\cR}{\mathcal R}

\newcommand{\cD}{\mathcal D}

\newcommand{\bbN}{\mathbb N}

\newcommand{\bbU}{\mathbb{U}}
\newcommand{\bbP}{\mathbb{P}}

\newcommand{\DN}{\Delta_N}

\newcommand{\dN}{\delta_N}
\newcommand{\dNt}{\delta_{N, t}}
\newcommand{\dNtp}{\delta_{N, t+1}}
\newcommand{\DM}{\Delta_M}
\newcommand{\hT}{h_\textrm{T}}
\newcommand{\JT}{J_\textrm{T}}
\newcommand{\dM}{\delta_M}

\newcommand{\blk}{\textrm{Blk}}
\newcommand{\Gte}{G_\text{tet}}
\newcommand{\Gted}{G_\text{ted}}
\newcommand{\qed}{\hfill\IEEEQEDhere}
\newcommand{\STGp}{STG$ ^+ $}
\newcommand{\deq}{\vcentcolon=}

\begin{document}	
\bstctlcite{IEEEexample:BSTcontrol}
%
\title{Finite-Horizon Optimal Control of Boolean Control Networks: A Unified Graph-Theoretical Approach}

\author{Shuhua~Gao,
	Changkai~Sun,
	Cheng~Xiang$ ^* $,
	Kairong~Qin,
	Tong Heng~Lee
	\thanks{Shuhua Gao, Tong Heng Lee, and Cheng Xiang (\textit{corresponding author}) are with the Department
		of Electrical \& Computer Engineering,  National University of Singapore, Singapore 119077, e-mail: elexc@nus.edu.sg.}
	\thanks{Kairong Qin is with the School of Optoelectronic Engineering and Instrumentation Science, and Changkai Sun is with the School of Biomedical Engineering, Dalian University of Technology, Dalian 116024, China.}
	}


%



\IEEEtitleabstractindextext{%
\begin{abstract}
This paper investigates the finite-horizon optimal control (FHOC) problem of Boolean control networks (BCNs) from a graph theory perspective. We first formulate two general problems to unify various special cases studied in the literature: (i) the horizon length is \textit{a priori} fixed; (ii) the horizon length is unspecified but finite for  given destination states. Notably, both problems can incorporate time-variant costs, which are rarely considered in existing work, and a variety of constraints. The existence of an optimal control sequence is analyzed under mild assumptions.  Motivated by BCNs' finite state space and control space, we approach the two general problems in an intuitive and efficient way under a graph-theoretical framework. A weighted state transition graph and its time-expanded variants are developed, and the equivalence between the FHOC problem and the shortest path problem in specific graphs is established rigorously. Two custom algorithms are developed to find the shortest path and construct the optimal control sequence with lower time complexity, though technically a classical shortest-path algorithm in graph theory is sufficient for all problems. Compared with existing algebraic methods, our graph-theoretical approach can achieve the state-of-the-art time efficiency while targeting the most general problems. Furthermore, our approach is the first one capable of solving Problem (ii) with time-variant costs. Finally, the Ara operon genetic network in \textit{E. coli} is used as a benchmark example to validate the effectiveness of our approach, and results of two tasks show that our approach can dramatically reduce the running time.
\end{abstract}

}

\maketitle

\IEEEdisplaynontitleabstractindextext

%
\IEEEpeerreviewmaketitle

\section{Introduction}  \label{sec: intro}
%
%
%

Boolean network (BN) was first proposed by Kauffman \cite{kauffman1969metabolic} to model gene regulatory networks, where each gene is assigned a Boolean variable to represent its expression state. The BN model has thereafter attracted increasing research interest in various fields, including studies on biomolecular networks in systems biology \cite{saadatpour2013boolean}, therapeutic interventions in clinical treatment \cite{datta2003external}, and the contagion dynamics during a financial crisis \cite{caetano2015boolean}, just to name a few.  In a BN, the binary variables interact with each other through Boolean functions, and exogenous (binary) inputs can be injected into these functions to affect the network dynamics, which is commonly referred to as a \textit{Boolean control network} (BCN) \cite{cheng2009controllability}.  In this study, we focus on finite-horizon optimal control (FHOC) of BCNs, which is useful for optimal therapeutic intervention strategy design in medical applications \cite{datta2003external}. Note that we view the \textit{finite horizon} in a more general sense: either a fixed horizon length, or an unknown but finite horizon length with respect to a given destination state. These two types of FHOC problems are referred to as \textit{fixed-time optimal control }and\textit{ fixed-destination optimal control} respectively in this paper.

In the last decade, a new matrix product called the \textit{semi-tensor product} (STP), which can convert a BN (BCN) into an algebraic state-space representation (ASSR), has been developed by Daizhan Cheng et al. \cite{cheng2010linear,zhao2010input}. The ASSR provides a systematic framework to study a wide range of control-theoretical problems related to BCNs, such as their controllability \cite{cheng2009controllability}, observability \cite{cheng2009controllability, laschov2013observability}, stabilization \cite{liang2017algorithms}, and various controller synthesis problems \cite{li2019robustness, zhao2015control}, among others. A number of optimal control problems for BCNs have been investigated in recent years using the STP and ASSR tools as well. In \cite{laschov2010maximum}, the Mayer-type optimal control problem (i.e., only terminal cost is considered) for single-input BCNs is addressed, and a necessary condition analogous to the Pontryagin's maximum principle is derived, which has been later extended to multi-input BCNs \cite{laschov2013pontryagin}. The minimum-energy control and minimum-time control of BCNs are investigated in \cite{li2013minimum} and \cite{laschov2013minimum} respectively. More general FHOC problems involving both stage cost and terminal cost are considered in \cite{fornasini2013optimal}, and the solution is given by a recursive algorithm as an analogy to the difference Riccati equation for discrete-time linear systems. More recently, Ref. \cite{zhu2018optimal}  targets the time-discounted stage cost and introduces a recursive algorithm based on a data structure called the optimal input-state transfer graph. The same problem is also investigated in \cite{cheng2015receding}, and a recursive solution for receding horizon optimal control of mix-valued probabilistic logical networks is obtained. In parallel to the study of FHOC problems, the more challenging infinite-horizon counterparts have also been attempted recently in several contributions using STP-based algebraic methods. For example, infinite-horizon optimal control with average cost is studied in \cite{zhao2010optimal, zhao2011floyd, fornasini2013optimal}, and  the time-discounted cost case is examined in \cite{cheng2014optimal, zhu2018optimal, wu2019optimal}.  

While we appreciate the above successes achieved with algebraic methods powered by the STP theory, one issue is that distinct methods are developed to solve different problems as reviewed above, even if we only consider FHOC. It appears that most existing work only deals with certain special cases of FHOC, for example, the minimum-energy control \cite{li2013minimum}, the minimum-time control \cite{laschov2013minimum}, the Mayer-type problem \cite{laschov2010maximum}, the two kinds of Lagrange-type problems \cite{cui2018optimal}, and the time-discounted finite-horizon problem \cite{cheng2015receding}, though they may share a lot in common. As mentioned in \cite{fornasini2013optimal} and \cite{zhang2017finite}, there are many other types of FHOC problems as well as lots of practical limitations, for instance, optimal control subject to various constraints \cite{faryabi2008optimal,zhang2017finite}, and more challenging problems with general time-variant costs. Consequently, a versatile approach instead of fragmented methods to solve all these common problems is highly desirable, which partly motivates this research. One goal of our study is to unify various tasks into an integrated framework and to develop systematic algorithms to handle the most general problems.

One serious concern about the algebraic approaches in BCN studies is the high computational complexity, since the number of operations grows exponentially with respect to the network size \cite{cheng2009controllability, zhao2010input, zhao2015control}. As a result, these approaches may quickly become computationally intractable as the number of state variables increase.  In fact, most control-theoretical problems of BCNs are NP-hard \cite{laschov2013minimum} due to the combinatorial nature of Boolean state variables, i.e., a BCN with $ n $ variables has up to $ N \deq 2^n  $ states, a phenomenon well-known as the \textit{curse of dimensionality}, such as NP-hardness of controllability \cite{akutsu2007control} and observability \cite{laschov2013observability}. However, this doesn't mean we are hopeless: many algebraic approaches run in a high-order polynomial time in terms of $ N $ (see Table \ref{tab: time complexity} for more details), which still leaves vast room for improvement by designing more efficient algorithms to decrease the order of the polynomial \cite{liang2017improved, liang2017algorithms}. This forms another objective of this research: we aim to reduce the computational complexity of FHOC for BCNs. 

To pursue more efficient algorithms, we notice that a BCN is characterized by its finite state space, finite control space, and deterministic state transitions. A direct consequence of this property is that the dynamics of a BCN can be adequately described by a state transition graph (STG). Going further, methods originating from graph theory appear to be promising for investigations of BCNs. For example,  computationally efficient methods for controllability \cite{liang2017improved, zhu2018further} and stabilization \cite{liang2017algorithms} of BCNs  all resort to certain graph-theoretical algorithms. 
Regarding the optimal control of BCNs, a couple of pioneering studies exist that attempt to combine the ASSR with tools from graph theory to improve computational efficiency. A typical example is the employment of Floyd-like algorithms, inspired by the Floyd-Warshall algorithm in graph theory, in both finite-horizon minimum-energy control \cite{li2013minimum} and infinite-horizon problems \cite{zhao2011floyd, cheng2014optimal}. An update-to-date study \cite{cui2018optimal} handles two kinds of Lagrange-type optimal control problems using the Dijkstra's algorithm instead. The above two algorithms are both initially developed to find shortest paths in a weighted graph. One limitation of both \cite{li2013minimum} and \cite{cui2018optimal} is that they only consider time-invariant costs and a special class of FHOC problems.
Motivated by these pioneering work borrowing tools from graph theory, we attempt to advance further and unify all common FHOC problems into an elegant and efficient graph-theoretical framework. 

The contributions of this paper are listed in three folds. First, we unify all common types of FHOC problems, including time-variant costs and various constraints, into two general problems depending on whether the horizon length is prespecified, which are subsequently reduced to shortest path problems by constructing particular variants of the STG. The correctness of this reduction is proved rigorously. Existence of an optimal control sequence is analyzed under mild conditions for both problems.
Second, we develop two intuitive algorithms to solve the above shortest path problems with superior efficiency. 
To be specific, only one algebraic method can achieve the same worst-case time complexity as ours, but our approach still tends to have better practical performance, which is demonstrated by two optimal control tasks for a genetic network in the bacteria \textit{E. coli}. Third, as far as we know, there are currently no published results on the fixed-destination optimal control problem with time-variant costs,  and our approach can handle it effectively using the identical graph-theoretical methodology. 

The remainder of this paper is organized as follows. In Section \ref{sec: preliminary}, we present some background knowledge about the STP, the ASSR, and the shortest path problem in graph theory. Section \ref{sec: problem formulation} introduces two general problems that can incorporate all specific FHOC problems studied in the current literature. In Section \ref{sec: stg}, the construction of the STG from a given initial state is discussed. After that, we detail the equivalence between the two FHOC problems and the shortest path problem in dedicated graphs and propose two efficient algorithms to solve them in Section \ref{sec: p1 solution} and \ref{sec: p2 solution} respectively. The time complexity of our approach is compared with that of existent work in Section \ref{sec: time complexity comparison}, and the practical running time of various methods to complete two tasks for the Ara operon network of \textit{E. coli} is measured in Section \ref{sec: benchmark}.  Finally, Section \ref{sec: conclusion} gives some concluding remarks.

\section{Preliminaries} \label{sec: preliminary}
\subsection{Notations}
For statement ease, the following notations  \cite{zhao2010input, cheng2010linear} are used.
\begin{enumerate}
	\item $ |\mathcal{S}| $ denotes the size (i.e., cardinality) of a set $ \mathcal{S} $.
	\item Let $\mathbb{R}$ and $\mathbb{N}$ denote the set of real numbers and nonnegative integers respectively. $ [l, r] \deq  \{l, l+1, \cdots, r - 1, r  \}$.
	\item $ f \ge B $ means a function $ f $ is bounded below by $ B $.
	\item $\mathcal{M}_{p\times q}$ denotes the set of all $p\times q$ matrices. Given $ A \in  $
	\item $ \textrm{Col}_i(M) $ denotes $ i $-th column of a matrix $ M $, and $ M_{ij} $ denotes the $ (i,j) $-th entry of $ M $.
	\item $ \delta_n^i \deq \textrm{Col}_i(I_n) $, where $ I_n \in \mathcal{M}_{n\times n} $ is the identity matrix. $\Delta_n \deq \{ \delta_n^i | i = 1, 2, \cdots, n \}$, and $ \Delta \deq \Delta_2 $. 
	\item A matrix  $M =[ \delta_{n}^{i_1}\; \delta_{n}^{i_2}\; \cdots\; \delta_{n}^{i_q} ] \in \mathcal{M}_{n\times q} $ with $\delta_n^{i_k}  \in \Delta_n, \forall k \in [1, q]$, is called a \textit{logical matrix}.  Let $\mathcal{L}_{n\times q} $ denote the set of all $ n\times q $ logical matrices.
	\item A matrix $ A \in  \mathcal{M}_{n\times mn}$ can be rewritten into a block form $ A = [\text{Blk}_1(A)\; \text{Blk}_2(A)\; \cdots\; \text{Blk}_m(A)] $, where $ \text{Blk}_i(A) \in \mathcal{M}_{n \times n}$ is the $ i $-th square block of $ A$.
	\item Logical operators \cite{cheng2010linear}: $\land$ for conjunction, $\lor$ for disjunction, $ \lnot $ for negation, and $\oplus $ for exclusive disjunction.
\end{enumerate}

\subsection{STP of Matrices and ASSR of BCNs}
This section revisits some necessary background knowledge about the STP and the ASSR of BCNs developed by Daizhan Cheng et al. \cite{cheng2010linear,zhao2010input,cheng2009controllability}. 
\begin{definition}\cite{zhao2010optimal}
	The STP of two matrices $ A \in   \mathcal{M}_{m\times n}$ and $ B \in   \mathcal{M}_{p\times q}$ is defined by
	\begin{equation*}
	A \ltimes B = (A \otimes I_{\frac{s}{n}})(B \otimes I_{\frac{s}{p}}),
	\end{equation*}
	where $\otimes$ denotes the Kronecker product,  and $ s $ is the least common multiple of $ n $ and $ p $. 
\end{definition}
\begin{remark}
	The STP is essentially a generalization of the standard matrix product, and all major properties of the standard matrix product remain valid under STP \cite{zhao2010optimal}. Thus, we will omit the symbol $\ltimes$ in the remainder when no confusion is caused. That is, all matrix products refer to STP by default.
\end{remark}

A logical function can be conveniently expressed in a multilinear form via STP. In this form, a Boolean value is identified by a vector as $ \text{TRUE} \sim \delta_2^1 $ and $ \text{FALSE} \sim \delta_2^2 $.

\begin{lemma} \cite{cheng2010linear} \label{lemma: structure matrix}
	Given any Boolean function $ f(x_1, x_2, \cdots, x_n): \Delta^n \rightarrow \Delta $, there exists a unique matrix $ M_f \in  \mathcal{L}_{2\times 2^n}$, called the \textit{structure matrix}, such that
	\begin{equation}
	f(x_1, x_2, \cdots, x_n) = M_f x_1 x_2 \cdots x_n.
	\end{equation}
\end{lemma}

Interested readers may refer to \cite{cheng2010linear} for computation of $ M_f $. For simplicity, set $ \ltimes_{i=1}^n A_i \deq A_1\ltimes\cdots\ltimes A_n $.

A general BCN $ \Sigma $ with $ n $ state variables (i.e., $ n $ nodes in the network) and $ m $ control inputs can be described as 
\begin{equation} \label{bcn}
\Sigma : 
\begin{cases}
x_1(t+1) = f_1(x_1(t), \cdots, x_n(t), u_1(t), \cdots, u_m(t))\\
x_2(t+1) = f_2(x_1(t), \cdots, x_n(t), u_1(t), \cdots, u_m(t))\\
\vdots \\
x_n(t+1) = f_n(x_1(t), \cdots, x_n(t), u_1(t), \cdots, u_m(t)),
\end{cases}
\end{equation}
where $ x_i(t) \in \Delta $ denotes the value of the $ i $-th variable at time $ t $, and $ f_i:  \Delta^{n+m} \rightarrow \Delta$ is the $ i $-th Boolean function, $ i\in [1, n]$. Additionally, $ u_j(t) \in \Delta $ denotes the $ j $-th control input at time $ t $, $ j \in [1, m]$. Set $ x(t) \deq \ltimes_{i=1}^n x_i(t)$ and $ u(t) \deq  \ltimes_{j=1}^m u_j(t) $. Clearly, we have $ x(t) \in \Delta_N $ and $ u(t) \in \Delta_M $, where $ N \deq 2^n $ and $ M \deq 2^m $.  $ N $ and $ M $ will be used through the whole text.

The ASSR of the BCN  in \eqref{bcn} is given by \cite{cheng2010linear,li2013minimum}:
\begin{equation} \label{eq: ASSR}
 \  x(t+1) = Lu(t)x(t),
\end{equation}
where $ L \in \mathcal{L}_{N \times MN} $, named the \textit{network transition matrix}, is computed by $ \text{Col}_j(L) = \ltimes_{i=1}^{n} \text{Col}_j(M_{f_i}), \forall j \in [1, MN]$, where $M_{f_i} \in  \mathcal{L}_{2\times MN}$ is the structure matrix of $ f_i $ in \eqref{bcn}. 

\subsection{Shortest Path Problem} \label{sec: SP}
The core of our graph-theoretical approach for FHOC is to transform the original problem into a shortest path problem in a particular graph and then locate the shortest path efficiently. We brief the shortest path problem in graph theory below.

Given a directed graph $ G=(V, E) $, where $ V = \{v_1, v_2, \cdots, v_n\} $ is a set of vertices, and $ E = \{(v_i, v_j) | v_i, v_j \in V \}$ is a set of directed edges, we assign each edge a real value, called its \textit{weight} or \textit{cost}.  Denote the weight of  the edge from $ v_i $ to $ v_j$ by $ w(v_i, v_j) $. Such a graph is known as a weighted directed graph. A path from a source vertex $ v_{i_0} \in V$ to a destination vertex $ v_{i_k} \in V $ is a sequence of vertices connected by edges, denoted by $ p=\left<v_{i_0}, v_{i_1},  \cdots, v_{i_k}\right> $. Let $ \epsilon(p) $ and $ |p| $ denote the number of edges and vertices in $ p $ respectively. If $ v_{i_0} = v_{i_k}$, $ p $ is a cycle. $ \epsilon(p) $ is also called the length of $ p $.

\begin{definition} \label{def: SP}
	The weight $ w(p) $ of a path $ p $  is the sum of the weights of its constituent edges:
	\begin{equation} \label{eq: path weight}
	w(p) = \sum_{j=0}^{k-1} w(v_{i_j}, v_{i_{j+1}}).
	\end{equation}
	A shortest path (SP) from $ v_{i_0} $ to $ v_{i_k} $ is any path  from $ v_{i_0} $ to $ v_{i_k} $ with the minimum weight among all possible paths.
\end{definition}

\section{Problem Formulation} \label{sec: problem formulation}
Despite the variety of FHOC problems of BCNs studied in the literature \cite{laschov2010maximum, laschov2013pontryagin, li2013minimum,laschov2013minimum,fornasini2013optimal,cheng2015receding,zhu2018optimal,cui2018optimal}, they can essentially be classified into two general types according to whether the horizon length is \textit{a priori }fixed or not.  This section will detail the mathematical formulation of both general problems. 

\subsection{Fixed-Time Optimal Control} \label{sec: problem 1}
Roughly speaking, our objective is to construct a fixed-length control sequence for the BCN \eqref{eq: ASSR} to optimize a given performance index \cite{faryabi2008optimal, fornasini2013optimal}. The scenario can become much more complicated in practice because of various constraints when designing control strategies. For example, we have to avoid dangerous states in therapeutic intervention \cite{datta2003external, faryabi2008optimal}. Besides, not all theoretical control inputs are practically realizable, e.g., we may still lack effective means to manipulate specific genes in a GRN, or a medical treatment is possibly unaffordable. Consequently, certain control inputs can be unavailable or just prohibited in specific states \cite{zhang2017finite}. Finally, it is common that we may want to steer the network to a particular state, e.g., to lead a gene regulatory network from a cancerous state (initial state) to a healthy state (terminal state) \cite{pal2006optimal, kim2013discovery}. Thus, we further enrich the fixed-time optimal control problem by constraining the terminal states.

A general problem reflecting the above intention subject to various constraints is formulated mathematically as follows.

\begin{problem} \label{problem: fixed length}
	Consider the BCN \eqref{eq: ASSR}, an initial state $ x_0 \in \Delta_N$, and a fixed horizon length $ T \in \mathbb{N}$. Fixed-time optimal control of a BCN is to determine an optimal control sequence of length $ T $ to the following optimization problem:
	\begin{align} \label{eq: general fixed length problem}
	\min_{u}\JT(u) = \hT(x(T)) + \sum_{t=0}^{T-1}g(x(t), u(t), t), \nonumber \\
	\textrm{s.t.} \begin{cases}
	x(t+1) = Lu(t)x(t) \\
	x(t) \in C_x \\
	u(t) \in C_u(x(t)) \\
	x(0) = x_0 \\
	x(T) \in \Omega
	\end{cases},
	\end{align}
	where $ u \deq (u(0), u(1), \cdots, u(T-1)) \in \DM^T$ is a control sequence; $ g: \Delta_N \times \Delta_M \times \bbN \rightarrow \mathbb{R} $ is the \textit{stage cost function}; and $ \hT:  \Delta_N \rightarrow \mathbb{R}$ is the \textit{terminal cost function}. $ C_x $ denotes the allowed states; $ C_u $ represents the state-dependent constraints on control inputs; and the terminal set $ \Omega \subseteq \DN $ denotes the set of desirable terminal states at time $ T $.
\end{problem}

\begin{remark}
	The horizon length $ T $ denotes a finite treatment window in therapeutic applications \cite{datta2003external, faryabi2008optimal}. Thus, optimal therapeutic intervention strategies can be developed within a treatment window by setting up a proper optimality criterion  $ \JT $ in \eqref{eq: general fixed length problem} by domain expert knowledge.
\end{remark}

\begin{remark}
	Most existing studies only deal with time-invariant stage cost, that is, the function $ g $ in \eqref{eq: general fixed length problem} doesn't really depend on time $ t $ (see \cite{zhao2011floyd, laschov2013minimum, cui2018optimal} for examples). Few studies consider a time-dependent $ g $ but only in restricted forms such as the time-discounted cost in \cite{zhu2018optimal} and \cite{cheng2015receding}. We intend to investigate the most general form \eqref{eq: general fixed length problem} directly, where the stage cost $ g $ can incorporate time $ t $ in any form.
\end{remark}

Problem \ref{problem: fixed length} represents a very general setting of fixed-time optimal control problems. A variety of specific finite-horizon problems investigated in existing work can be rewritten easily into this general form. We demonstrate in the sequel how to specialize Problem \ref{problem: fixed length} to some specific fixed-time optimal control problems with different characteristics in the literature.
\begin{enumerate}
	\item \textit{No terminal state constraints} \cite{fornasini2013optimal, zhu2018optimal, cui2018optimal}. If the terminal state is not specified, we only need $ \Omega = \Delta_N $.
	\item \textit{A single desired terminal state} \cite{laschov2013minimum, li2013minimum}. Set $ \Omega = \{ x_d \} $, where $ x_d $ is the single desired terminal state.  
	\item \textit{Only stage cost} \cite{li2013minimum, cheng2015receding,cui2018optimal}. To remove the terminal cost, simply nullify it by $ \hT(\cdot) \equiv 0$.
	\item \textit{Only terminal cost (also known as the Mayer problem)} \cite{laschov2010maximum, laschov2013pontryagin}. Just set zero stage costs: $ g(\cdot, \cdot, \cdot) \equiv 0$.
	\item \textit{Special terminal or stage cost functions} \cite{fornasini2013optimal, li2013minimum, cheng2015receding}.  For example, \cite{zhu2018optimal} and \cite{cheng2015receding} consider the time-discounted finite-horizon optimal control, which can be expressed by \eqref{eq: general fixed length problem} with $ g(x(t), u(t), t) = \lambda^t c_g(u(t), x(t)) $, where $ 0 < \lambda < 1$ is the discount factor, and $ c_g(\cdot, \cdot) $ is a time-invariant stage cost. As another example, the energy function in \cite{li2013minimum} is obtained through $ g(x(t), u(t), t) = u^\top(t)Qu(t) $, and the more general quadratic cost function in \cite{fornasini2013optimal} can be implemented straightforwardly by
	\begin{align*}
	&\hT(x(T)) = x(T)^\top Q_h x(T), \\
	&g(x(t), u(t), t) = \begin{bmatrix} x(t)^\top u(t)^\top \end{bmatrix}
	\begin{bmatrix} 
	Q & S \\
	S^\top & R 
	\end{bmatrix} \begin{bmatrix} x(t) \\ u(t) \end{bmatrix},
	\end{align*}
	where $ Q_h,  Q, S, \text{and } R $ are proper weight matrices. 
\end{enumerate}

\subsection{Fixed-Destination Optimal Control} \label{sec: fixed-dest problem}
A common task that Problem \ref{problem: fixed length} cannot cover is time-optimal control, which aims to find a control sequence to drive the BCN from a given initial state $ x_0 $ to another destination state $ x_d $ in minimum time. Such time-optimal control has been widely studied for traditional linear time-invariant (LTI) discrete-time systems, such as the famous deadbeat controller (see \cite{o1981discrete} for a review). In \cite{laschov2013minimum}, the minimum-time control of BCNs is first investigated, and the same problem is considered again in  \cite{chen2016minimum} but with impulsive disturbances. Such optimality concept can be generalized to other criteria, for example, the minimum-energy control studied in \cite{li2013minimum}, which attempts to steer the BCN to a target state using minimum energy.

Note that although the horizon length is not  \textit{a priori} fixed, it must be finite for a well-posed problem, that is, the specified destination state $ x_d $ should be reachable from the initial state $ x_0 $ in finite steps. It is the main reason we consider such problems as another class of FHOC problems in a general sense. To further generalize this problem, just like Problem \ref{problem: fixed length}, we make it admit a set of destination states and subject to various constraints, which is formalized as follows. 
\begin{problem}\label{problem: unknown length}
	Consider the BCN \eqref{eq: ASSR}, an initial state $ x_0 \in \Delta_N$, and a terminal set $ \Omega \subseteq \Delta_N$. The fixed-destination optimal control problem of a BCN is to determine an optimal control sequence of a variable length to the optimization problem,
	\begin{align} \label{eq: general unknown length problem}
	&\min_{u} J(u) = h(x(K), K) + \sum_{t=0}^{K-1} g(x(t), u(t), t), \nonumber\\
	&\textrm{s.t.} \begin{cases}
	x(t+1) = Lu(t)x(t) \\
	x(t) \in C_x \\
	u(t) \in C_u(x(t)) \\
	x(0) = x_0 \\
	x(K) \in \Omega
	\end{cases},
	\end{align}
	where $ u = (u(0), u(1), \cdots, u(K-1)) \in \DM^K$ represents a control sequence, and $ K \in \mathbb{N}$ indicates an unknown but finite horizon length.  The terminal cost function $ h: \Delta_N \times\mathbb{N} \rightarrow \mathbb{R} $ and the stage cost function $g: \Delta_N \times \Delta_M \times \bbN \rightarrow \mathbb{R} $ can be time-dependent. $ C_x $, $ C_u $, and $ \Omega $ denote the state constraint, the control constraint, and desirable destination states respectively.
\end{problem}

\begin{remark}
	Despite the apparent similarity in mathematical forms between Problem  \ref{problem: fixed length} and Problem \ref{problem: unknown length}, they will be treated by distinctly different means to maximize the computational efficiency of each problem. In existing work like \cite{li2013minimum} and \cite{cui2018optimal}, only specific problems with a time-independent stage cost function and a single destination state are investigated, and no terminal cost is considered. However, we argue that it is sensible to set different costs if the desired terminal state is reached at different time, for example, when the desired state refers to a good state of physical health.
\end{remark}

Compared with the fixed-time optimal control (Problem \ref{problem: fixed length}), there are fewer studies on Problem \ref{problem: unknown length}. As far as we know, only the following two custom problems have been investigated in the literature, both with a single destination state, i.e., $ \Omega=\{ x_d \} $, and a  time-invariant stage cost function. 
\begin{enumerate}
	\item \textit{Minimum-time control} \cite{laschov2013minimum, chen2016minimum}. This kind of control is often referred to as time-optimal control. In Problem \ref{problem: unknown length}, simply set $ g(\cdot, \cdot, \cdot) \equiv 1 $ and $ h (\cdot, \cdot) = 0$, and we have $ J(u) = K$. Thus, we are indeed minimizing the number of steps to steer the BCN \eqref{eq: ASSR} from $ x_0 $ to $ x_d $.
	\item \textit{Minimum-energy control} \cite{li2013minimum}. Let $ g(x(t), u(t), t) = u^\top(t) Q u(t) $ and $ h (\cdot, \cdot) = 0$, where $ Q  $ is a positive definite diagonal matrix measuring energy consumption. 
\end{enumerate}

\section{State Transition Graph} \label{sec: stg}
In this section,we introduce the state transition graph (STG) of a BCN, which is useful for both Problem \ref{problem: fixed length} and \ref{problem: unknown length}.  An efficient algorithm based on bread-first search (BFS) \cite{cormen2009introduction} in a graph is developed to construct the STG.

\subsection{Reachability of a BCN}
To construct the STG, we need to first discuss the reachability of a BCN \cite{cheng2009controllability,zhao2010input}, especially the reachable set of a given intial state, whose definition is given as follows.

\begin{definition} \label{def: reachable set}
	The set of states that can be reached from $x(0) = x_0 \in \Delta_N$ at  time $ t= d $ is $ \mathcal{R}_d(x_0) $. Given a set $ X \subseteq \DN $, let $ \mathcal{R}_d(X) \deq \cup_{x \in X} \mathcal{R}_d(x) $ for notational simplicity. The complete reachable set of a given state $ x_0$ is $ \cR(x_0) \deq \cup_{d \in \bbN}\cR_d(x_0) $.
\end{definition}

$ \cR(x_0) $ can be obtained algebraically by calculating powers of the network transition matrix $ L $ \cite{cheng2009controllability, zhao2010input}. However,  from a graph-theoretical view, a computationally economical way is to adopt the standard bread-first search (BFS)  procedure to build $ \cR(x_0) $ iteratively \cite{cormen2009introduction,liang2017algorithms}. More interestingly, this method is closely related with the adjacency-list representation of a graph \cite{cormen2009introduction}, and only requires successive computation of the one-step reachable set $ \cR_1(\cdot) $. 
Recall that $ u(t) \in \DM $ and $ x(t)  \in \DN$ are both logical vectors with all zero entries except a single entry of value 1. Thus, given $ u(t) = \dM^k $, we have $ Lu(t) = \blk_k(L) $. The following lemma follows immediately.

\begin{lemma} \label{lemma: R1}
	The one-step reachable set of a state $ \dN^i $ for a BCN \eqref{eq: ASSR} with no constraints is
	\begin{equation} \label{eq: R1}
	\mathcal{R}_1(\delta_N^i) = \{  \textrm{Col}_i(\textrm{Blk}_k(L)) |  k \in [1, M]  \},
	\end{equation}
	and, if there are constraints like Problem \ref{problem: fixed length} (or Problem \ref{problem: unknown length}), 
	\begin{equation} \label{eq: R1-constrained}
	\mathcal{R}_1(\delta_N^i) = \{  \textrm{Col}_i(\textrm{Blk}_k(L)) |  k \in [1, M], \dM^k \in C_u(\dN^i)  \} \cap C_x.
	\end{equation}
\end{lemma}

Note that we ignore the terminal constraint in \eqref{eq: R1-constrained} and will handle it later when building specific graphs. It is possible that more than one control input can attain the transition from $ \dN^i $ to $ \dN^j $. Collect these qualified control inputs into a set $ U^{ij} $:
\begin{equation}  \label{eq: U^ij}
	U^{ij} = \{ \dM^k \in C_u(\dN^i)| k \in [1, M], \textrm{Col}_i(\textrm{Blk}_k(L)) = \dN^j  \}.
\end{equation}

By Definition \ref{def: reachable set} and Lemma \ref{lemma: R1}, the $ d- $step reachable set of a state $ \dN^i $ is obtained efficiently with the recursion:
\begin{equation}
\cR_d(\dN^i) = \cR_1(\cR_{d-1}(\dN^i)), d \ge 2. 
\end{equation}

\subsection{Construction of an STG}
As aforementioned, a BCN is characterized by its finite state space, where the finite control inputs coordinate the state transitions deterministically. This feature makes it possible to capture the complete dynamics of a BCN with a directed graph, called the \textit{state transition graph} (STG). 
\begin{definition} \label{def: STG}
	Consider the BCN \eqref{eq: ASSR}.
	Its state transition graph (STG) is a directed graph $ G =(V, E) $, where $ V = \Delta_N $ is the vertex set, i.e., one vertex for each state, and the edge set is
	\begin{equation}  \label{eq: edges of STG}
	E = \{(\delta_N^i, \delta_N^j) | \delta_N^i \in V, \dN^j \in \cR_1(\delta_N^i)   \},
	\end{equation}
	i.e., one edge for each one-step transition between states.
\end{definition}

In this study, we only care about states reachable from an initial state $ x_0 $. We denote such an STG by $ G =(V, E, x_0) $ with $ V = \mathcal{R}(x_0)$. Hereafter we will use the terms \textit{vertex} and \textit{state} interchangeably when no ambiguity is caused.

In addition, we may assign a weight to each edge of the STG corresponding to the cost of the each state transition, which will be detailed in following sections.

\begin{remark}
	The \textit{input-state graph} (in a matrix form), proposed by \cite{cheng2014optimal, zhao2010input}, uses an input-state pair $ (\delta_M^i, \delta_N^j) $ as an vertex, leading to $ MN$ vertices in total. Our STG with only $ N $ or $\cR(x_0) $ vertices is potentially more space and time efficient. 
\end{remark}

Note from \eqref{eq: edges of STG} that $ \cR_1(\dN^i) $ essentially denotes the successors of the vertex $ \dN^i $ in the STG. We thus can build the STG following a BFS procedure and get $ \cR(x_0) $ at the same time. BFS is a graph traversal algorithm where the neighbors of a vertex are visited in a FIFO (first-in-first-out) order. Algorithm \ref{alg: STG} details the construction of an STG. Running Algorithm \ref{alg: STG}, we get the reachable set $ \cR(x_0) $, i.e., $ R $ in the algorithm, and the adjacency list of each vertex, $ \cR_1(x), \forall x \in \cR(x_0)$,  which effectively gives the STG $ G =(V, E, x_0) $ by Definition \ref{def: STG}.

\begin{algorithm}[tb]
	\caption{Construction of State Transition Graph via BFS} \label{alg: STG}
	\begin{algorithmic}[1]
		\Input a BCN \eqref{eq: ASSR} and relevant constraints if any
		\Output Adjacency-list representation of an STG
		\State Initialize a FIFO queue $ Q $ and a set $ R $
		\State Append $ x_0 $ to $ Q $ and $ R $
		\While{$ Q $ is not empty}
			\State $ \dN^i \gets $ $\textrm{dequeue}(Q) $
			\ForAll{$ \dN^j \in  \cR_1(\dN^i)$} \label{line: for 1}  \Comment See Lemma \ref{lemma: R1}
				\If{$ \dN^j \notin R $}
					\State Append $ \dN^j $ to $ Q $ and $ R $
				\EndIf
			\EndFor
		\EndWhile
	\end{algorithmic}
\end{algorithm}

\textit{Time complexity analysis of Algorithm \ref{alg: STG}.} For each state $ \dN^i \in \cR(x_0)$, Eq. \eqref{eq: R1} or \eqref{eq: R1-constrained} needs $ M $ operations, and obviously $ |\cR_1(\dN^i)| \le M $. The \textbf{for} loop has at most $ M $ iterations; and the \textbf{while} loop runs $ |V| $ times because each vertex is enqueued and dequeued exactly once. Thus, the time complexity of Algorithm \ref{alg: STG} is $ O(M|V|) $, which is equivalent to $ O(MN) $ since there are at most $ N $ states (vertices). This conforms to the celebrated theorem that BFS runs in linear time with respect to the number of edges and vertices \cite{cormen2009introduction}, i.e., $ O(|V|+|E|) $, because the STG has at most $ MN $ edges. 

\begin{remark}
	Algorithm \ref{alg: STG} implies that, any state $ x  \in \cR(x_0)$ can be reached from $ x_0 $ in less than $ N $ steps, i.e., $ \cR(x_0) = \cup_{0 \le d < N} \cR_d(x_0) $. This result is intuitive: if a trajectory from $ x_0 $ to $ x $ contains more than $ N $ states, there must be repetitive ones, and we can remove such cycles to shorten the trajectory.
\end{remark}

\begin{example} \label{example: 1}
	Consider a BCN with $ m = 2 $ control inputs and $ n = 3 $ state variables \cite{li2013minimum} as follows:
	\begin{equation}	\label{eq: example 1 bcn}
		\Sigma_1:
		\begin{cases}
		x_1(t+1) = x_2(t) \land (u_1(t) \oplus x_3(t) ) \\
		x_2(t+1) = \lnot x_1(t) \\
		x_3(t+1) = u_2(t) \oplus x_2(t)
		\end{cases}.
	\end{equation}
	Its ASSR \eqref{eq: ASSR} has a $ 8 \times 32 $ transition matrix $ L $, which is omitted here to conserve space. For illustration purpose, we set up the state constraints and control constraints arbitrarily below:
	
	\begin{equation} \label{eq: example 1 constraints}
		\begin{cases}
			C_x &= \Delta_8 \setminus  \{\delta_8^8\} \\
			C_u(x) &= \begin{cases}
			\{ \delta_4^3, \delta_4^4  \}, \  \textrm{if } x = \delta_8^6 \\
			\Delta_4 \setminus \{\delta_4^2\}, \ \textrm{otherwise}
			\end{cases}
		\end{cases}
	\end{equation}
	i.e., the state $ \delta_8^8 $ is forbidden; the control $ \delta_{4}^2 $ is unavailable to all states; and only control $ \{ \delta_4^3, \delta_4^4  \} $ is applicable to state $ \delta_8^6 $.
\end{example}

\begin{figure}[tb]
	\centering
	\includegraphics[width=35mm]{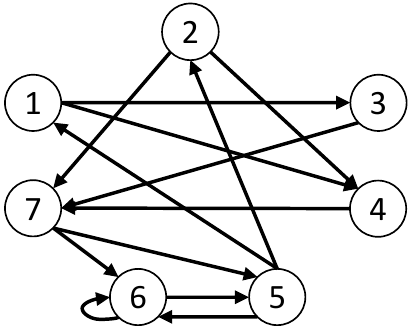} 
	\caption{ STG of the BCN \eqref{eq: example 1 bcn} in Example \ref{example: 1}. Each circle labeled $ i $ represents a state $ \delta_{8}^i $, and the arrows denote state transitions.}
	\label{fig: STG of example 1}
\end{figure}

The STG $ G=(V, E, x_0) $ of $ \Sigma_1 $ with states reachable from $ x_0 = \delta_8^1 $ subject to the given constraints \eqref{eq: example 1 constraints} is obtained by Algorithm \ref{alg: STG} and shown in Fig. \ref{fig: STG of example 1}.  \qed

\section{Solve Problem \ref{problem: fixed length} via Dynamic Programming} \label{sec: p1 solution}
In this section, we first validate the existence of optimal solutions to Problem \ref{problem: fixed length}, and then we introduce a time-expanded variant of the STG and develop an efficient algorithm via dynamic programming to solve Problem \ref{problem: fixed length} based on that graph.

\subsection{Existence of Optimal Solutions}

As we will discuss next, the stage costs and the terminal costs in \eqref{eq: general fixed length problem} are mapped to edge weights of a graph to reduce Problem \ref{problem: fixed length} to an SP problem. It is known that the SP problem is well defined only if the graph contains no negative-weight cycles \cite{cormen2009introduction}. Nonetheless, since the length of a path in Problem \ref{problem: fixed length} is fixed due to the fixed horizon, this condition is no longer required. Only the following assumption is needed for Problem \ref{problem: fixed length}  to avoid an  optimal value of negative infinity.

\begin{assumption} \label{ass: bounded below}
	The stage cost function $ g $ and the terminal cost function $ \hT $ in Problem \ref{problem: fixed length} are both bounded from below.
\end{assumption}

Now we are ready to present the following conclusion about the existence of optimal solutions to Problem \ref{problem: fixed length}.

\begin{proposition} \label{prop: 1}
	Consider Problem \ref{problem: fixed length} under Assumption \ref{ass: bounded below}. An optimal control sequence $ u^* $ to \eqref{eq: general fixed length problem} exists if and only if $ \Omega \cap \mathcal{R}_T(x_0) \ne \emptyset$.
\end{proposition}

\begin{IEEEproof}
	(\textit{Necessity})
	$ \mathcal{R}_T(x_0) $ includes all states reachable from $ x_0 $ at $ t = T $  subject to the constraints in \eqref{eq: general fixed length problem} (except $ \Omega $), i.e., $ x(T) \in  \mathcal{R}_T(x_0), \forall u \in \Delta_M^T$. If $ \Omega \cap \mathcal{R}_T(x_0) = \emptyset$, no feasible solutions exist for Problem \ref{problem: fixed length}, i.e., Problem \ref{problem: fixed length} is infeasible. Thus, the necessity of $ \Omega \cap \mathcal{R}_T(x_0) \ne \emptyset$ is obvious.
	
	(\textit{Sufficiency}) Note that the solution space $ \mathbb{U} $ of Problem \ref{problem: fixed length} is of finite size, which contains at most $ M^T $ candidate solutions. Besides, $ \Omega \cap \mathcal{R}_T(x_0) \ne \emptyset$ implies that at least one feasible solution $ u \in  \mathbb{U}$ exists which can steer the BCN from $ x_0 $ to a terminal state $ x_f \in \Omega $ at time $ T $ under constraints. Moreover, Assumption \ref{ass: bounded below} ensures that $ \JT(u) $ is bounded from below, $ \forall u \in \mathbb{U} $.
	A straightforward exhaustive enumeration of $ \mathbb{U}$  can yield the optimal solution $ u^* $ satisfying  $ \JT(u^*) = \min_{u \in \mathbb{U}} \JT(u) $. 
\end{IEEEproof}

\begin{remark}
	If there are no constraints in Problem \ref{problem: fixed length}, i.e., $C_x \equiv \DN, C_u(\cdot)  \equiv \DM, \Omega \equiv \Delta_N $, an optimal control sequence always exists, which is widely studied in existing work, e.g.,  \cite{cui2018optimal, zhu2018optimal}. Note additionally that the optimal control sequence $ u^* $ may not be unique.
\end{remark}

\subsection{Time-Expanded Fixed-Time State Transition Graph}  \label{sec: data structure problem 1}
Fig. \ref{fig: STG of example 1} shows that each edge in the STG corresponds to a state transition of the BCN, whose weight indicates the transition cost. This fact motivates us to connect Problem \ref{problem: fixed length} to the shortest-path (SP) problem on the STG. However, in contrast to the standard SP problem in graph theory,  Problem \ref{problem: fixed length} poses three substantial challenges. First, the number of time steps is fixed to $ T $, that is, we want only $ T $-edge paths. Second, the stage cost function $ g $ is time-dependent, indicating that the edge weights may vary with time.  Finally, there is an additional terminal cost given by $ \hT $. Consequently, the classic SP algorithms can no longer be applied. To overcome these obstacles, we get inspiration from the space-time network in dynamic transportation network studies  \cite{pallottino1998shortest} and propose a new graph called the \textit{Time-Expanded fixed-Time State Transition Graph} (TET-STG), which attaches timestamp to state transitions by stretching the STG along the time dimension. Besides, we introduce a \textit{pseudo-state} $ \dN^0 $ to handle terminal states and their costs. A formal definition is given below.

\begin{definition}  \label{def: TET-STG}
	Consider Problem \ref{problem: fixed length}. The TET-STG $ \Gte = (V, E, x_0, T) $ is a weighted directed graph constructed by:
	\begin{itemize}
		\item $ V = \cup_{t=0}^{T+1} V_t $, where $ V_t = \cR_t(x_0), \forall  t \in [0, T-1] $, $ V_{T} = \cR_T(x_0) \cap \Omega $, and $ V_{T+1} = \{ \dN^0  \} $.
		\item $ E = \cup_{t=0}^{T} E_t$, where $ E_t = \{(\delta_N^i, \delta_N^j) | \delta_N^i \in V_t,  \delta_N^j \in V_{t+1} \cap \mathcal{R}_1(\delta_N^i)    \},  \forall t \in [0, T-1]$, and \\ $ E_T = \{ (\delta_N^i, \dN^0) | \dN^i \in V_T \} $.
	\end{itemize}
	Note that a state may appear at multiple time points, but they are treated as distinctive
	vertices from a graph perspective. Denote the vertex representing the state $ \dN^i $ at time $ t $ by $ \dNt^i $. The weight of the edge $ (\dNt^i,  \delta_{N, t+1}^{j}) \in E_t $ is 
	\begin{equation} \label{eq: w_t^ij}
	w(\delta_{N, t}^{i}, \delta_{N, t+1}^{j}) = 
		\begin{cases}
		 \min_{\dM^k \in U^{ij}} g(\delta_N^{i}, \dM^k, t), t \in [0, T - 1]  \\
		 \hT(\delta_N^i), t = T
		\end{cases}
	\end{equation}
	where $ U^{ij} $ is given in \eqref{eq: U^ij}, and the unique vertex at time $ T+ 1$ refers to the pseudo-state $ \delta_{N, T+1}^0 $. Denote the control that achieves the weight (cost) in \eqref{eq: w_t^ij} by $ u_t^{ij} $:
	\begin{equation} \label{eq: u_t^ij}
		u^{ij}_t = \argmin_{\dM^k \in U^{ij}} g(\delta_N^{i}, \dM^k, t), t \in [0, T - 1].
	\end{equation}
\end{definition}

\begin{remark}
	We note a slight abuse of notations in the above: $ \cR_t(x_0) $ includes time information implicitly, and thus $ V_i \cap V_j = \emptyset, \forall i \ne j $.
	Besides, recall Eq. \eqref{eq: U^ij}: among the possibly nonunique control inputs that enable a state transition, we pick definitely the one of  lowest cost in \eqref{eq: w_t^ij} for optimal control purpose. The role of the pseudo-state $ \delta_{N, T+1}^0 $ is to incorporate the terminal cost into the graph at the \textit{pseudo-time} $ T+1 $.
\end{remark}

Despite its seemingly complex definition, the TET-STG can be built handily by acquiring $V_0, E_0, V_1, \cdots, E_{T-1}, V_T$ successively similar to the BFS in Algorithm \ref{alg: STG}. In practical implementation, the one-step transition between states need to computed only once: supposing there is a transition $ (\dN^i, \dN^j)  $ in the STG, i.e., $ \dN^j \in \cR_1(\dN^i) $,  if we have a vertex $ \dNt^i \in V_t $, then there exists a succeeding vertex $ \dNt^j \in V_{t+1} $ and an edge $ (\dNt^i, \dNtp^j) \in E_t$ in the TET-STG. 

\begin{figure}[tb]
	\centering
	\includegraphics[width=75mm]{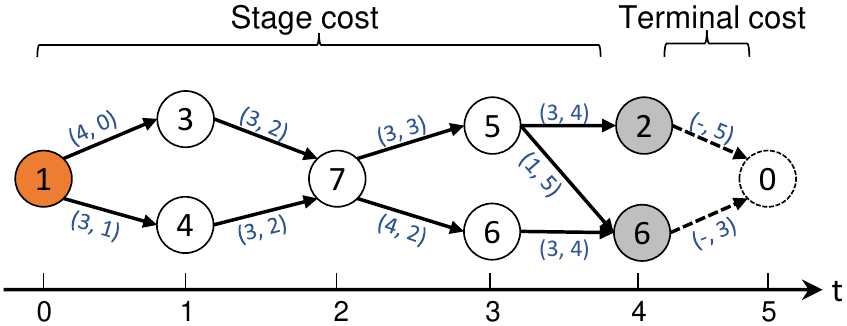} 
	\caption{ TET-STG in Example \ref{example: 2} with the initial state $ x_0 = \delta_8^1$ (the orange circle). A circle labeled $ i $ at time $ t $ denotes a vertex $ \delta_{8, t}^i $. We highlight the terminal states in $ \Omega $ by gray circles and the pseudo-state $ \delta_8^0 $ by a dashed circle and dashed edges. The annotation of each edge $ (k, w) $ means that this transition is achieved with control $ \delta_4^k $ at a cost of $ w $ (i.e., the edge's weight). Note that the dashed edges going into the pseudo-state need no control.}
	\label{fig: TET-STG of example 2}
\end{figure}

\begin{example} \label{example: 2}
	Consider the BCN $ \Sigma_1 $ \eqref{eq: example 1 bcn} in Example \ref{example: 1} again. In addition to the constraints \eqref{eq: example 1 constraints}, Problem \ref{problem: fixed length} is set up by $T = 4,  x_0 = \{ \delta_8^1 \} $, and $ \Omega = \{ \delta_8^2, \delta_8^6 \} $. The costs are:
	\begin{equation*}
		g(x(t), u(t), t) = u(t)^{\top}Qu(t) + t, \ h_\text{T}(x(T)) = x(T)^{\top}Rx(T),
	\end{equation*}
	where $ Q = \textrm{diag}(2, 3, 1, 0) $ and $ R = \textrm{diag}(3, 5, 4, 0, 1, 3, 6, 0) $. 
	
	Take the transition $ (\delta_{8}^3, \delta_{8}^7) $ as an example. We have $ U^{37} = \{ \delta_4^1, \delta_{4}^3  \} $, which justifies the use of \eqref{eq: u_t^ij}. The obtained TET-STG is shown in Fig. \ref{fig: TET-STG of example 2}. Now it is clear that, though the stage cost function $ g $ itself is time-dependent, the weight of each edge in the TET-STG becomes time-invariant after we expand the STG along the time axis.
	\qed
\end{example}

Recall Proposition \ref{prop: 1}. From the construction of the TET-STG, it is clear that an optimal control sequence exists if $ V_T \ne \emptyset$, because each path from $ \delta_{N, 0}^{i_0} $ to $\delta_{N, T+1}^{0} $ yields a feasible solution to Problem \ref{problem: fixed length}, where $ \delta_{N}^{i_0}  = x_0 $. 

As aforementioned, we will transform Problem \ref{problem: fixed length} into an SP problem in the TET-STG. First, note that each one-step transition in the TET-STG is already attained with a minimum cost thanks to \eqref{eq: w_t^ij} and \eqref{eq: u_t^ij}. A direct consequence is the following lemma whose correctness is intuitive.

\begin{lemma} \label{lemma: trajectory & path}
	Consider Problem \ref{problem: fixed length} under Assumption \ref{ass: bounded below}. For any feasible control sequence $ u = (\dM^{k_0}, \dM^{k_1}, \cdots, \dM^{k_{T-1}} ) $ that steers the BCN along a state trajectory $ s = (\dN^{i_0}, \dN^{i_1}, \cdots, \dN^{i_T} ) $ with $x_0 = \delta_{N}^{i_0}  $ and $ \dN^{i_T} \in \Omega $, there holds $ \JT(u)  \ge w(p) $, where $ p $ is the corresponding path in the associated TET-STG,
	\begin{equation} \label{eq: a path}
		p = \big< \delta_{N, 0}^{i_0}, \delta_{N, 1}^{i_1}, \cdots, \delta_{N, T}^{i_T}, \delta_{N, T+1}^{0}  \big>,
	\end{equation}
	and $ \JT(u)  = w(p) $ holds if $ \dM^{k_t} = u_t^{i_t i_{t+1}}, \forall t \in [0, T - 1]$.
\end{lemma}

\begin{IEEEproof}
	According to edge weights in \eqref{eq: w_t^ij} and the optimal one-step control in \eqref{eq: u_t^ij}, the weight of the path $ p $ in \eqref{eq: a path} is 
	\begin{equation*}
		w(p) = \hT(\delta_N^{i_T}) + \sum_{t=0}^{T-1} g(\delta_N^{i_t}, u_t^{i_t i_{t+1}}, t).
	\end{equation*}
	From the definition of $ \JT$ in \eqref{eq: general fixed length problem}, we have 
	\begin{equation*}
		\JT(u) = \hT(\delta_N^{i_T}) + \sum_{t=0}^{T-1}g(\delta_N^{i_t}, \dM^{k_t}, t).
	\end{equation*}
	We can see that the correctness of Lemma \ref{lemma: trajectory & path} is obvious.
\end{IEEEproof}

The following theorem establishes the connection between fixed-time optimal control and the SP in the TET-STG.
\begin{theorem} \label{thm: problem 1 SP}
	Consider Problem \ref{problem: fixed length} under Assumption \ref{ass: bounded below} and suppose it is feasible. Given $ x_0 = \delta_{N}^{i_0}$, if an SP from $ \delta_{N, 0}^{i_0} $ to $\delta_{N, T+1}^{0} $ in the TET-STG $ \Gte = (V, E, x_0, T) $ is 
	\begin{equation} \label{eq: p^* of p1}
		p^* =  \big< \delta_{N, 0}^{i_0}, \delta_{N, 1}^{i_1^*}, \cdots, \delta_{N, T}^{i_T^*}, \delta_{N, T+1}^{0}  \big>,
	\end{equation}
	then the minimum value of the cost function \eqref{eq: general fixed length problem} is $ \JT^* = w(p^*)$, and an optimal control sequence is 
	\begin{equation} \label{eq: U^* of p1}
	u^* = \{ u_0^{i_0 i_1^*}, u_1^{i_1^*i_2^*}, \cdots, u_{T-1}^{i_{T-1}^* i_T^*}    \}.
	\end{equation}
\end{theorem}

\begin{IEEEproof}
	Suppose the solution space of Problem \ref{problem: fixed length} is $ \bbU \ne \emptyset$, and the set of paths from $ \delta_{N, 0}^{i_0} $ to $\delta_{N, T+1}^{0} $ in $ \Gte $ is $ \bbP $. There holds $ \bbP \ne \emptyset $ because each state trajectory driven by $ u \in \bbU $ corresponds to a path in $ \bbP $. Moreover, we have $ |\bbP| \le |\bbU| $ since different control sequences may lead to identical state trajectories (see \eqref{eq: U^ij}). As we have shown in the proof of Proposition \ref{prop: p2}, $ |\bbU| $ is finite, which implies $ |\bbP| $ is also finite. Additionally, Assumption \ref{ass: bounded below} and Eq. \eqref{eq: w_t^ij} ensures that all edge weights in $ \Gte $ are bounded from below, which guarantees $ w(p) $ is bounded from below for any $ p \in \bbP $ because $ p $ has exactly $ T+1 $ edges.  Thus, there must exist a shortest path $ p^* $.
	
	Given any $ u \in \bbU $, Lemma \ref{lemma: trajectory & path} tells that there exists $ p \in \bbP $ such that $ \JT(u) \ge w(p) \ge w(p^*) $. Furthermore, for the control sequence $ u^* \in \bbU $  in \eqref{eq: U^* of p1} that attains $ p^* $, it leads to a state trajectory $ s^* = (\delta_{N}^{i_0}, \delta_{N}^{i_1^*}, \cdots, \delta_{N}^{i_T^*}) $ with $ \dN^{i_T^*} \in \Omega $ according to the TET-STG in Definition \ref{def: TET-STG}. Thus, $ u^* $ is a feasible solution to Problem \ref{problem: fixed length}. Furthermore,
	Lemma \ref{lemma: trajectory & path}  states $ \JT(u^*) = w(p^*)$. Hence, we have $ \JT(u) \ge \JT(u^*) $, $\forall u \in \bbU$.
\end{IEEEproof}

\subsection{Dynamic Programming (DP) in TET-STG}
The problem  following Theorem \ref{thm: problem 1 SP} immediately is how to locate an SP from $\delta_{N, 0}^{i_0}  $ to $ \delta_{N, T+1}^{0} $ in the TET-STG efficiently (see \eqref{eq: p^* of p1}). Since the layered structure of the TET-STG doesn't contain any cycles, the classic SP algorithms \cite{cormen2009introduction} such as the Floyd-Warshall algorithm and the Dijkstra's algorithm can be applied directly. However, in view of the fixed horizon length in Problem \ref{problem: fixed length}, we propose a custom method based on dynamic programming (DP) to achieve better time efficiency, which can even beat the state-of-the-art Dijkstra's algorithm (Remark \ref{remark: algorithm 1}).

The intuition behind our DP approach is that any sub-path of an SP is itself an SP as well \cite{cormen2009introduction}. Such optimal
sub-structure is a strong indicator that DP based methods are applicable. The following theorem formalizes this idea. 

\begin{theorem} \label{thm: DP}
	Consider Problem \ref{problem: fixed length} under Assumption \ref{ass: bounded below} and suppose it is feasible. In its associated TET-STG $ \Gte = (V, E, x_0, T) $ with $x_0 = \delta_N^{i_0}  $, let $ F(\delta_{N, t}^j) $ denote the weight of an SP from vertex $ \delta_{N, 0}^{i_0} $ to vertex $ \delta_{N, t}^{j} $, $ \forall t \in [0, T + 1] $, and let $ P(\delta_{N, t}^j) $ denote the predecessors of $ \delta_{N, t}^j $ in $ \Gte $:
	\begin{equation} \label{eq: predecessors}
	P(\delta_{N, t}^j) = \{ \delta_{N, t-1}^i | (\delta_{N, t-1}^i, \delta_{N, t}^j) \in E  \}, t \in [1, T + 1]
	\end{equation}
	Then $ F(\delta_{N, t}^j) $ can be obtained by the following recursion:
	\begin{equation} \label{eq: DP recursion}
	F(\delta_{N, t}^j) = \min_{\delta_{N, t-1}^i   \in P(\delta_{N, t}^j) } F(\delta_{N, t-1}^i ) + w(\delta_{N, t-1}^i, \delta_{N, t}^j),
	\end{equation}
	for $ t \in [1, T + 1] $, 
	and the base condition is $ F(\delta_{N, 0}^{i_0}) = 0 $. If $ p^* $ is an SP from $ \delta_{N, 0}^{i_0} $ to $ \delta_{N, T+1}^{0} $, we have $ w(p^*) = F(\delta^0_{N, T+1}) $.
\end{theorem}

\begin{IEEEproof}
	We first show the correctness of the recursion \eqref{eq: DP recursion} by induction. Due to the layered structure of $ \Gte $, any vertex $ \delta_{N, t}^j \in V$ can only be reached from a certain vertex $ \delta_{N, t-1}^i \in V $ in one step, $ t \ge 1 $. Assume $ F(\delta_{N, t-1}^i ), t\ge 1, $ represents the minimum weight from vertex $ \delta_{N, 0}^{i_0} $ to vertex $ \delta_{N, t-1}^i $. For any  path from  $ \delta_{N, 0}^{i_0} $ to $ \delta_{N, t}^j $ that passes $ \delta_{N, t-1}^i $, its minimum weight is obviously $ F(\delta_{N, t-1}^i ) $ + $ w(\delta_{N, t-1}^i, \delta_{N, t}^j) $. Eq. \eqref{eq: DP recursion} examines all such  predecessors of $ \delta_{N, t}^j $, and the minimum value $ F(\delta_{N, t}^j) $  is clearly the minimum weight of any path from  $ \delta_{N, 0}^{i_0} $ to $ \delta_{N, t}^j $. Besides, for $ t=0 $, the base case $ F(\delta_{N, 0}^{i_0}) = 0 $ is clearly true. Thus, we have verified that  $F(\delta^j_{N, t}) $ is the minimum weight of any path from $ \delta_{N, 0}^{i_0} $ to $ \delta^j_{N, t} $. It implies directly that $F(\delta^0_{N, T+1}) $ is the weight of the SP $ p* $ from $ \delta_{N, 0}^{i_0} $ to $ \delta_{N, T+1}^{0} $.
\end{IEEEproof}

\begin{remark}
	Eq. \eqref{eq: DP recursion} is essentially a form of the Bellman optimality equation, a widely applied tool in optimal control \cite{datta2003external}, though we state it from a graph-theoretical perspective.
\end{remark}

Combing Theorem \ref{thm: problem 1 SP} and \ref{thm: DP}, we get the minimum cost $ \JT^* $ for Problem \ref{problem: fixed length} by $ \JT^* = w(p^*) =  F(\delta^0_{N, T+1})$. Nonetheless, we are more interested in the optimal control sequence $ u^* $ that attains $ \JT^* $. The key is to record the minimizer to \eqref{eq: DP recursion}: if $ \dNt^j $ is in an SP, then the minimizer $ \delta_{N, t-1}^{i^*} $ must be in the SP as well. At last, $ p^* $ can be reconstructed accordingly, followed by $ u^* $ acquired with \eqref{eq: U^* of p1}. 
Algorithm \ref{alg: fixed horizon} implements the DP method by Theorem \ref{thm: DP} to solve $ F(\delta_{N, t}^j) $ and to reconstruct the optimal control sequence. Note that, for maximal time efficiency, the predecessors \eqref{eq: predecessors} of each vertex are stored while building $ \Gte $. 

\begin{algorithm}[tb]
	\caption{Fixed-Time Optimal Control via DP} \label{alg: fixed horizon}
	
	\begin{algorithmic}[1] 
		\Input $ L, T, \hT, g, C_x, C_u, \Omega $, and $ x_0 $ in Problem $ \ref{problem: fixed length} $
		\Output The optimal control sequence $ u^* $ and $ \JT^* $
		\State Build the TET-STG $ \Gte$ according to Definition \ref{def: TET-STG}

		\Function{ShortestPath}{$t, j, \cM, F$}
		\If {$ t=0 $}  \Comment Base condition
		\State \Return 0
		\EndIf
		\If {$ F $ has the key $ (t, j) $}        \Comment Memoization
		\State \Return $ F[(t, j)] $
		\EndIf
		\State $ i^* \leftarrow 0$, $\  d^* \leftarrow \infty $ 
		\ForAll {$ \delta_{N, t-1}^i \in P(\delta_{N, t}^j) $}   \Comment Recursion \eqref{eq: DP recursion}
		\State
		\begin{varwidth}[t]{\linewidth}
						$d \leftarrow  $ \textproc{ShortestPath}($ t - 1, i , \cM, F $) + \par
						\hskip\algorithmicindent \hspace{2pt} $ w(\delta_{N, t-1}^i, \delta_{N, t}^j) $ 
				\end{varwidth}
		\If {$ d< d^* $}
		\State $ i^* \leftarrow i $, $\quad  d^*  \leftarrow d$
		\EndIf
		\EndFor
		\State $ F[(t, j)] \leftarrow d^* $
		\State $ \cM[(t, j)] \gets i^*$	\Comment Record the minimizer of \eqref{eq: DP recursion}
		\State \Return $ d^* $
		\EndFunction
		
		\LineComment  Call the recursive function \textproc{ShortestPath}
		\State Initialize two dictionaries (i.e., hash tables) $ \cM $ and $ F $
		\State $ \JT^* \gets$ \Call{ShortestPath}{$T + 1, 0 , \cM, F$}  \label{line: SP}
		\State Create an array $ u^* $ of length $ T $
		\State $ j \gets \cM[(T+1, 0)], \quad t \gets T$
		\While{$ t> 0 $}
			\State $ i \gets \cM[(t, j)], \quad u^*[t - 1] \gets u_{t-1}^{ij} $	\Comment See \eqref{eq: u_t^ij}
			\State $ j \gets i, \quad t \gets t - 1 $
		\EndWhile
		
	\end{algorithmic}
\end{algorithm}

\textit{Time complexity analysis of Algorithm \ref{alg: fixed horizon}.} $ \Gte $ has $ T + 2 $ layers with at most $ Z \deq |\cR(x_0)| $ vertices per layer and at most $ MZ $ edges between any two successive layers. Like Algorithm \ref{alg: STG}, $ \Gte $ is constructed via BFS in linear time $ O(TMZ) $. The core function of Algorithm \ref{alg: fixed horizon}, \Call{ShortestPath}{}, implements top-down DP via the memoization technique \cite{cormen2009introduction}, i.e., storing results into $ F $ and retrieving the cached results if same inputs recur. Memoization can avoid repetitive computation \cite{cormen2009introduction}, and accordingly each edge of $ \Gte$ is processed only once. \Call{ShortestPath}{} in Line \ref{line: SP} thus takes time $ O(TMZ) $. The remaining construction of $ u^* $ runs in $ O(T) $. The overall worst-case time complexity of Algorithm \ref{alg: fixed horizon} is hence $ O(TMZ) $, or equivalently, $ O(TMN) $, because we always have $ Z \le N $.

\begin{remark} \label{remark: algorithm 1}
	The recent work \cite{cui2018optimal} establishes a weighted directed graph to formulate the $ k $-edge shortest path problem as well. However, neither the time-dependent stage cost function nor the terminal constraint set is considered in \cite{cui2018optimal}. Moreover, the work \cite{cui2018optimal} directly applies the standard Dijkstra's algorithm, which is designed for general SP problems and less efficient than our Algorithm \ref{alg: fixed horizon} in solving Problem \ref{problem: fixed length}, whose running time is $ O(TMN + TN\log(TN)) $ instead.
\end{remark}

\begin{example} \label{example: p1}
	Recall Example \ref{example: 2}. Running Algorithm \ref{alg: fixed horizon}, we get the minimum value of the cost function $ \JT^* = 11 $, and the optimal control sequence $ u^* = (\delta_4^4, \delta_4^3, \delta_4^4, \delta_4^3) $. The correctness of this result can be easily verified by inspecting Fig. \ref{fig: TET-STG of example 2} and enumerating all paths from $ \delta_{4, 0}^1 $ to $ \delta_{4, 5}^0 $. \qed
\end{example}

\section{Solve Problem \ref{problem: unknown length} via Dijkstra's Algorithm} \label{sec: p2 solution}
In this section, the existence of optimal solutions to Problem \ref{problem: unknown length} is first examined.  Then, we divide Problem \ref{problem: unknown length} into two cases depending on whether cost functions are time-dependent.  Both cases will be conquered by Dijkstra's algorithm, but different graph structures are constructed to maximize efficiency.

\subsection{Existence of Optimal Solutions}
Like Problem \ref{problem: fixed length}, we use the terminal costs and the stage costs of Problem \ref{problem: unknown length} as the edge weights of specific graphs. However, the major difference is that the number of state transitions in Problem \ref{problem: unknown length} is not fixed. Consequently, the condition that no negative-weight cycles exist in any state trajectory from $ x_0 $ to $ x_d \in \Omega$ is mandatory \cite{cormen2009introduction}; otherwise, the cost $ J $ can always be reduced by repeating a negative-weight cycle, and no SP exists. We thus require the following conditions to guarantee existence of optimal solutions to Problem \ref{problem: unknown length}.

\begin{assumption} \label{ass: nonnegative}
	The cost functions in Problem \ref{problem: unknown length} satisfy three conditions: (i) $ h $ is bounded from below;  (ii) $ g $ is nonnegative; (iii) $ h $ and $ g $ are both nondecreasing with respect to time $ t $, i.e., $ h(\dN^i, t_2) \ge h(\dN^i, t_1),  \forall t_2 > t_1, \forall \dN^i \in \Omega$, and $ g(\dN^i, \dM^k, t_2) \ge g(\dN^i, \dM^k, t_1), \forall t_2 > t_1, \forall (\dN^i, \dM^k) \in \Delta_N \times \DM$.
\end{assumption}

The rationality of condition (i) is obvious, just like Problem \ref{problem: fixed length}, to ensure a finite optimal value. Condition (ii) can be technically relaxed to the nonexistence of negative-weight cycles, though it would be quite difficult to verify such a condition in practice because edge weights (i.e., stage costs) can vary with time.
We can justify condition (iii) intuitively by imaging a special scenario. Suppose a state trajectory from $ x_0  $ to $ x_d \in \Omega$ contains a cycle of zero weight. Then a possible result is that the more cycling the BCN  does along this cycle, the more the cost criterion $ J $ can be reduced, once $ g $ or $ h $ can decrease as time passes. Consequently, the optimal control sequence does not have a finite length. Note that if $ g $ and $ h $ do not depend on $ t $, which is the most common case in the literature, condition (iii) is satisfied naturally. 

\begin{remark} \label{rmk: h nonngeative}
	Following \cite{zhu2018optimal, fornasini2013optimal}, we can always assume $ h$ is  nonnegative without affecting the optimal solution. Supposing $ h $ is bounded from below by a constant $ B_h $, the optimal control sequence to \eqref{eq: general unknown length problem} is the same one that minimizes $ J'(u)  = J(u) - B_h$, and we have a nonnegative terminal cost function now: $ h'(x(K), K) = h(x(K), K) - B_h \ge 0$.
\end{remark}

The following proposition confirms the existence of an optimal solution to Problem \ref{problem: unknown length} under the above conditions.
\begin{proposition} \label{prop: p2}
	Consider Problem \ref{problem: unknown length} under Assumption \ref{ass: nonnegative}.  There exists an optimal control sequence $ u^* $ satisfying $ |u^*| < |\cR(x_0)| $ that minimizes \eqref{eq: general unknown length problem} if and only if $ \Omega \cap \mathcal{R}(x_0) \ne \emptyset$.
\end{proposition}

\begin{IEEEproof}
	(\textit{Necessity)} Since $ \mathcal{R}(x_0) $ denotes all states reachable from $ x_0 $, no state $ x_d \in \Omega $ can be reached if $ \Omega \cap \mathcal{R}(x_0) = \emptyset$.
	
	(\textit{Sufficiency}) If $ \Omega \cap \mathcal{R}(x_0) \ne \emptyset$, there must exist control sequences that steer the BCN from $ x_0 $ to a state $ x_d \in \Omega$. Suppose such a control sequence is $ u =\big(u(t)\big)_{t=0}^{k-1} $, and the resultant trajectory is $ s =\big(x(t)\big)_{t=0}^{k}  $, where  $x(0) = x_0 $ and $ x(k) = x_d$. Furthermore, we claim that if $|u| = k \ge  |\mathcal{R}(x_0)| $,  there must exist a shorter control sequence $ \bar{u} $ such that $ |\bar{u}| < |\mathcal{R}(x_0)|  $ and $ J(\bar{u})  \le J(u)$. This claim is proved below.
	
	$|u| \ge  |\mathcal{R}(x_0)| $ implies that $ |s| > |\mathcal{R}(x_0)| $, which means that $ s$ contains repetitive states because $ x \in  \mathcal{R}(x_0)$ for any state $ x $ in $ s $.  Assume one such pair of repetitive states is $ x(i) = x(j)$, $ 0 \le i < j \le k$, and thus $c = ( x(i), x(i+1), \cdots, x(j) )  $ is  a circular sub-trajectory. We can remove this cycle (except $ x(j) $) from $ s $, and obviously the remaining states $ s' = (x(0), \cdots, x(i-1), x(j), x(j+1), \cdots, x(k-1), x(k))$  still constitute a trajectory from $ x_0 $ to $ x_d $ driven by a shortened control sequence $ u' = (u(0), \cdots, u(i-1), u(j), u(j+1), \cdots, u(k-1)) $. Let $ r \deq j -i  > 0$, and it follows that
	\begin{align} \label{eq: JU}
	J(u) &-J(u') = h(x_d, k) - h(x_d, k - r) + \sum_{t=i}^{j - 1} g(x(t), u(t), t)  \nonumber\\
	&+ \sum_{t=j}^{k - 1} \left[ g(x(t), u(t), t) - g\big(x(t), u(t), t - r\big) \right].                                    
	\end{align}
	
	Condition (ii) and (iii) in Assumption \ref{ass: nonnegative} guarantee  \eqref{eq: JU} is nonnegative, i.e., $J (u') \le  J(u) $. The above cycle elimination procedure can be repeated until a control sequence $ \bar{u} $ satisfying $ |\bar{u}| < |\mathcal{R}(x_0)| $ is obtained. We have certainly $ J(\bar{u}) \le J(u)$.
	
	The above claim implies that it is enough to consider the candidate set $ \bar{\mathbb{U}} = \{u| |u| < |\cR(x_0)|, u\in \mathbb{U}\} $ for an optimal solution, where $ \mathbb{U}$ is the feasible set of Problem \ref{problem: unknown length}. Obviously, the set $ \bar{\mathbb{U}} $ is finite, and Assumption \ref{ass: nonnegative} guarantees $ J(u) $ is bounded from below, $ \forall u \in \bar{\mathbb{U}} $. Therefore,  an optimal solution $ u^* \in \bar{\mathbb{U}}$ to Problem \ref{problem: unknown length} must exist such that $ |u^*| < |\cR(x_0)|$. 
\end{IEEEproof}

\subsection{Case 1: Time-Invariant Stage Cost and Terminal Cost} \label{sec: p2 case1} 
In this case, neither $ g $ nor $ h $ of \eqref{eq: general unknown length problem} depends on time $ t $: the STG becomes a static graph, whose edge weights are permanently fixed. Furthermore, if there is only one destination state with zero terminal cost, Problem \ref{problem: unknown length} degrades to a standard SP problem in the STG. This simplest case has been solved in \cite{li2013minimum, laschov2013minimum, cui2018optimal}. We address the more general problems here, where multiple destination states with non-zero terminal costs are allowed. Following the same idea in solving Problem \ref{problem: fixed length}, we introduce an extra pseudo-state $ \delta_{N}^0 $ as well as the terminal set into the STG and term the new graph \STGp.

\begin{definition} \label{def: STG+}
	Consider Case 1 of Problem \ref{problem: unknown length}. The extended state transition graph STG$^+ $, denoted by $ G^+ = (V, E, x_0) $, is an extension of the STG constructed as follows:
	\begin{enumerate}
		\item Build the STG $ G = (V, E, x_0) $ by Definition \ref{def: STG}, and the weights are assigned like \eqref{eq: w_t^ij}, though time-independent:
		\begin{equation*}
			w(\delta_N^{i}, \delta_N^{j}) = \min_{\dM^k \in U^{ij}} g(\delta_N^{i}, \dM^k, \alpha), \ (\delta_N^{i}, \delta_N^{j})  \in E.
		\end{equation*}
		\item Add into $ G $ the pseudo-state $ \dN^0 $ and its incoming edges,
		\begin{equation*}
			V \gets V \cup \{ \dN^0 \}, \quad E \gets E \cup E^0
		\end{equation*}
		where $ E^0 = \{ (\delta_N^i, \delta_N^0) | \delta_N^i \in V \cap \Omega  \} $ with weights: 
		\begin{equation*}
			w(\delta_N^i, \delta_N^0) = h(\delta_N^i, \alpha), \ (\delta_N^i, \delta_N^0)  \in E^0
		\end{equation*}
		Here $ \alpha \in \mathbb{N}$ is an arbitrary integer used as a placeholder.
	\end{enumerate}
\end{definition}

Since all edge weights of $ G^+ $ are time-invariant, we denote the weight and the associated control input of each edge $ (\delta_N^i, \delta_N^j) \in E $ by $ w^{ij} $ and $ u^{ij}  $ respectively. Of course, the incoming edges of $ \dN^0 $ needs on control input. 

\begin{figure}[tb]
	\centering
	\includegraphics[width=45mm]{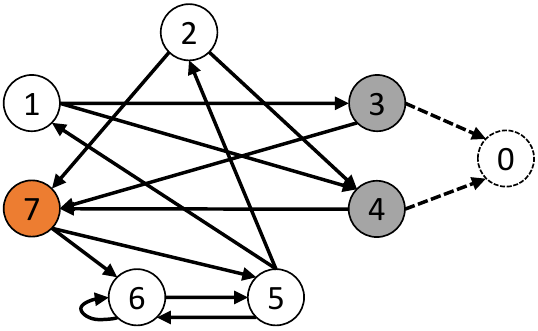} 
	\caption{ \STGp in Example \ref{example: STG+} with the initial state $ x_0 = \delta_8^7 $ (in orange) and two destination states $ \{ \delta_{8}^3, \delta_{8}^4 \} $ (in gray). The weight and the control associated with each edge are not shown here for clarity. Note that the dashed edges going into the pseudo-state $ \delta_{8}^0 $ need no control.}
	\label{fig: example STG+}
\end{figure}

\begin{example} \label{example: STG+}
	Consider Example \ref{example: 1} again: the BCN \eqref{eq: example 1 bcn} is subject to the constraints \eqref{eq: example 1 constraints}. Suppose the desired terminal states are $ \Omega = \{ \delta_{8}^3, \delta_{8}^4 \} $, and the initial state is $ x_0 = \delta_{8}^7 $. We assign the following time-invariant cost:
	\begin{equation*}
		g(x(t), u(t), t) = x(t)^{\top}Qx(t) + u(t)^{\top}Ru(t) , \ h(x(K), K) \equiv 0,
	\end{equation*}
	where $ Q = \textrm{diag}(2, 5, 1, 4, 1, 3, 6, 0) $ and $ R = \textrm{diag}(0, 3, 1, 4) $.
	The \STGp of this case can be easily built on basis of the STG in Fig. \ref{fig: STG of example 1} following Definition \ref{def: STG+}, which is shown in Fig. \ref{fig: example STG+}. As we see, the \STGp   is akin to the STG but with an additional pseudo-state and edge weights assigned.
	\qed
\end{example}

To transform Problem \ref{problem: unknown length} in this case into an SP problem in the \STGp, we first give the following lemma.
\begin{lemma} \label{lemma: path p2c1}
	Consider Case 1 of Problem \ref{problem: unknown length} under Assumption \ref{ass: nonnegative} with time-independent $ g $ and $ h $. Given any control sequence $ u = (\dM^{k_0}, \dM^{k_1}, \cdots, \dM^{k_{K-1}} ) $ of an unknown length $ K $ that steers the BCN along a state trajectory $ s = (\dN^{i_0}, \dN^{i_1}, \cdots, \dN^{i_K} ) $ with $x_0 = \delta_{N}^{i_0}  $ and $ \dN^{i_K} \in \Omega$, there holds $ J(u)  \ge w(p) $, where $ p $ is the corresponding path in the associated \STGp,
	\begin{equation} \label{eq: a path STG+}
	p = \big< \delta_{N}^{i_0}, \delta_{N}^{i_1}, \cdots, \delta_{N}^{i_K}, \delta_{N}^{0}  \big>,
	\end{equation}
	and $ J(u)  = w(p) $ holds if $ \dM^{k_t} = u^{i_t i_{t+1}}, \forall t \in [0, K - 1]$.
\end{lemma}
\begin{IEEEproof}
	With the edge weights given by Definition \ref{def: STG+}, the weight of the path $ p $ is
	\begin{equation*}
		w(p) = h(\dN^{i_K}, K) + \sum_{t=0}^{K - 1} g(\dN^{i_t}, u^{i_t i_{t+1}}, t).
	\end{equation*}
	Note that $ u^{i_t i_{t+1}} $ is the optimal control input to transit the BCN from $ \dN^{i_t} $ to $ \dN^{i_{t+1}} $ in one step, i.e., $g(\dN^{i_t}, u^{i_t i_{t+1}}, t) \le g(\dN^{i_t}, \dM^{k_t}, t), \forall t \in [0, K - 1]  $. Comparing $ w(p) $ with the cost $ J(u) $ in \eqref{eq: general unknown length problem}, Lemma \ref{lemma: path p2c1} is obviously true.
\end{IEEEproof}

The following theorem relates Case 1 of Problem \ref{problem: unknown length} to the SP problem in an \STGp.
\begin{theorem}\label{thm: p2 case 1}
	Consider Case 1 of Problem \ref{problem: unknown length} under Assumption \ref{ass: nonnegative} with time-independent $ g $ and $ h $, and suppose it is feasible. In the STG$ ^+ $ $ G^+ = (V, E, x_0), x_0=\delta_N^{i_0} $,  there exists an SP $ p^* = \big< \delta_N^{i_0}, \delta_N^{i_1^*},  \delta_N^{i_2^*}, \cdots, \delta_N^{i_{K-1}^*}, \delta_N^{i_K^*}, \delta_N^0 \big>, K \le |V| - 2, $   such that the minimum value of the cost function \eqref{eq: general unknown length problem} is $ J^* = w(p^*) $. The corresponding optimal control sequence is $ u^* = \{ u^{i_0i_1^*},  u^{i_1^*i_2^*}, \cdots, u^{i_{K-1}^* i_K^*} \} $.
\end{theorem}

\begin{IEEEproof}
	Since the problem is feasible, there exists control sequences that steers the BCN from the initial state $ x_0 $ to a destination state $ x_d \in \Omega$, implying that there exist paths from $ \dN^{i_0} $ to $ \dN^0 $ in $ G^+ $. Note that $ G^+ $ has no negative cycles because (i) Assumption \ref{ass: nonnegative} states $ g $ is nonnegative and (ii) the pseudo-state $ \dN^0 $ has only incoming edges, whose weights assigned by $ h $ may be negative though. With this fact, we can easily show that an SP $ p^* $ from $ \dN^{i_0} $ to $ \dN^0 $ with at most $ |V| $ vertices exists following the cycle elimination procedure in Proposition \ref{prop: p2}. In fact, this is a fundamental theorem in graph theory \cite{cormen2009introduction}. 
	
	Recall Lemma \ref{lemma: path p2c1}: given any feasible control sequence $ u $, there exists a path $ p $ from $ \dN^{i_0} $ to $ \dN^0 $ in $ G^+ $ such that $ J(u) \ge w(p) \ge w(p^*) = J(u^*)$. The above theorem holds clearly.
\end{IEEEproof}

Theorem \ref{thm: p2 case 1} has reduced Problem \ref{problem: unknown length} with time-invariant costs to a regular \textit{single-pair} SP problem \cite{cormen2009introduction} from $ x_0 $ to $\delta_N^0 $ in the \STGp.  Since $ g $ is nonnegative and the pseudo-state $ \dN^0 $ has only incoming edges, 
We have discussed in Remark \ref{rmk: h nonngeative} that we can always assume $ h $ is nonnegative without loss of generality, though it is not mandatory. Since $ g $ is also nonnegative, the \STGp  has only nonnegative edge weights. The fastest known SP algorithm for such graphs is Dijkstra's algorithm \cite{cormen2009introduction}.

We make two modifications to the normal implementation of Dijkstra's algorithm for this optimal control problem. First, like Algorithm \ref{alg: fixed horizon}, the vertices that compose the SP are recorded to reconstruct the optimal control sequence later. Second, we terminate the search process once vertex $\delta_N^0$ is reached because we are only interested in the SP from $ x_0 $ to $\delta_N^0$. Algorithm~\ref{alg: dijkstra} presents the modified Dijkstra's algorithm to solve Case 1 of Problem \ref{problem: unknown length}. Since Dijkstra's algorithm is a well-established SP algorithm in graph theory (see \cite[Chapter~24]{cormen2009introduction} for details), the proof of its correctness is omitted here.

\begin{algorithm}
	\caption{Fixed-destination optimal control with time-invariant costs using modified Dijkstra's algorithm}
	\label{alg: dijkstra}
	\begin{algorithmic}[1]
		\Input $ L, h, g, C_x, C_u, \Omega $, and $ x_0 = \dN^{i_0}$ in Problem \ref{problem: unknown length}
		\Output The optimal control sequence $ u^* $ and $ J^* $
		
		\State Build the STG$ ^+ $ $ G^+=(V, E, x_0) $ by Definition \ref{def: STG+} \label{alg2line: stg}
		\State Create a min-priority queue $ \mathcal{Q} $ and two dictionaries $ \cM,  \cD $  
		
		\ForAll{$\delta_N^i \in  V$ }  \Comment Set initial distances
		\IfThenElse{$ i = i_0 $} {$ \cD [i] \leftarrow 0$}{$ \cD [i] \leftarrow \infty$}
		\State Add $ i $ into $ \mathcal{Q} $ with its priority $ \cD [i] $
		\EndFor

		\While{$ |\mathcal{Q}| > 0 $}  \Comment Continue until we reach $ \dN^0 $ \label{line: main start}
		\State $ i \leftarrow$ extract the minimum-priority item from $ \mathcal{Q} $
		\IfThen{$ i = 0 $}{\textbf{break}}  \Comment Early termination
		
		\ForAll{$ j \in \{ j' | (\dN^i, \dN^{j'}) \in E  \}$}  \label{alg2line: neighbor}
		\State $ d \leftarrow \cD[i] + w^{ij}$
		\If {$ d < \cD[j] $}
		\State $\cD[j] \leftarrow d, \  \cM[j] \leftarrow i$
		\State Update the priority of $ j $ in $ \mathcal{Q} $ to $ d $
		\EndIf
		\EndFor
		\EndWhile  \label{line: main end}
		\State $ J^* \gets \cD[0] $		\Comment Minimum weight from $ \dN^{i_0} $ to $ \dN^0 $   \label{alg3line: dijkstra end}
 		\State Create an empty array $ u^* $, and  $ j \leftarrow 0 $ \label{line: reconstruct}
		\While{$ j \ne i_0 $}
		\State $ i \leftarrow \cM[j] $  \Comment Edge $ (\delta_N^i, \delta_N^j) $ is in the SP
		\State Append $ u^{ij} $ \eqref{eq: u_t^ij} to $ u^* $    except $ j = 0 $    
		\State $ j \leftarrow i $
		\EndWhile
		\State Reverse $ u^* $ in place \label{line: reconstruct end}
	\end{algorithmic}
\end{algorithm}

\textit{Time complexity analysis of Algorithm \ref{alg: dijkstra}.} A key data structure in Dijkstra's algorithm is the priority queue \cite{cormen2009introduction}, in which each item has a \textit{priority} and the one with highest (or lowest) priority is first served. If the priority queue is implemented with a Fibonacci heap, Dijkstra's algorithm has a running time of $ O(|E| + |V|\log|V|) $ \cite{cormen2009introduction}. In the beginning, just like the STG, the construction of the \STGp takes time $ O(M|V|) $. The more complicated Dijkstra's part (Line 2 to \ref{alg3line: dijkstra end}) runs in $ O(M|V|+ |V|\log |V|) $ accordingly. Finally, since the SP found by Dijkstra's algorithm contains at most $ |V| $ vertices, the construction of $ u^* $ (Line \ref{line: reconstruct} to \ref{line: reconstruct end}) runs in $ O(|V|) $. Overall,  the worst-case time complexity of Algorithm \ref{alg: dijkstra} is $ O(M|V| + |V|\log |V|) $, or equivalently, $ O(MN + N\log N) $.

\begin{example}
	We test Algorithm \ref{alg: dijkstra} with Example \ref{example: STG+}. The optimal value  is $ J^* = 13 $ and the optimal control sequence is $ u^* = ( \delta_4^1, \delta_4^3, \delta_4^1 )$. The state trajectory of the BCN is thus $ s^* = ( \delta_{8}^7, \delta_{8}^5, \delta_{8}^2, \delta_{8}^4 ) $. It is easy to verify in Fig. \ref{fig: example STG+} that an SP from $ \delta_{8}^7 $ to  $ \delta_{8}^0 $ is  $ p^* = \big<  \delta_{8}^7, \delta_{8}^5, \delta_{8}^2, \delta_{8}^4,  \delta_{8}^0 \big> $ and $ w(p^*) = 13 $.
	\qed
\end{example}

\subsection{Case 2: Time-Variant Stage Cost and Terminal Cost} \label{sec: p2 case2}
As aforementioned, the classic SP algorithms will not work once the edge weights may vary with time. To the best of knowledge, there are still no published studies on fixed-destination optimal control of BCNs with time-varying costs. Recall the TET-STG proposed in Section \ref{sec: p1 solution}, and we naturally attempt to handle this time-variant case for Problem \ref{problem: unknown length} in a similar way. However, one immediate difficulty is that, unlike Problem \ref{problem: fixed length}, the horizon length is not known beforehand, which prevents the reuse of the DP-based Algorithm \ref{alg: fixed horizon} directly.

Hopefully, we may resort to Proposition \ref{prop: p2} to overcome this obstacle: it is sufficient to consider only control sequences of size less than $ |\cR(x_0) |$ to find the optimal one, though the exact length remains unknown. On the other hand, Algorithm \ref{alg: fixed horizon} works once a horizon length $ T $ is given. 
A straightforward workaround that reuses Algorithm \ref{alg: fixed horizon} to solve this case thus comes to our mind as follows.
\begin{procedure} \label{proc: case 2 every n}
	Solve Case 2  by Reusing Algorithm \ref{alg: fixed horizon}.
	\begin{itemize}
		\item Step 1. Given $ T \in [1, |\cR(x_0)| - 1] $, transform this case into Problem \ref{problem: fixed length} by setting $ \hT(x(T))  = h(x(T), T) $.
		\item Step 2. Solve the problem obtained above with Algorithm \ref{alg: fixed horizon} to get an optimal control sequence $ u_T^*, |u_T^*| = T $. 
		\item Step 3. Repeat Step 1 and Step 2 for all possible $ T $'s, and finally yield $ u^* = \argmin_{u_T^*} (J(u_T^*)) $.
	\end{itemize}
\end{procedure}

The enumeration of all possible horizon lengths like Procedure \ref{proc: case 2 every n} is essentially the idea underlying the algebraic approach \cite{li2013minimum} for minimum-energy control towards a given target state, though it only considers time-invariant stage cost. A similar idea is also adopted in \cite{zhao2011floyd} to detect the minimum average-weight cycle in the input-state space. However, such a somewhat brute-force method is still inevitably computationally expensive, e.g., the running time is $ O(N^4) $ in \cite{li2013minimum, zhao2011floyd}. In our case, even though each subproblem for a specific $ T $ can be solved by the more efficient Algorithm \ref{alg: fixed horizon}, the overall time complexity of Procedure \ref{proc: case 2 every n} is still as high as $ O(MN^3) $. 

Since we only care about control sequences that has a size less than $ |\cR(x_0)| $, to further reduce the computational burden, we devise another approach by adapting the TET-STG to an unknown but limited horizon length. More interestingly, Algorithm \ref{alg: dijkstra} initially developed for Case 1 can be reused on the resultant graph. We call this new graph a \textit{Time-Expanded fixed-Destination State Transition Graph} (TED-STG). It construction is similar to the TET-STG detailed as follows.

\begin{definition} \label{def: TED-STG}
	Consider Case 2 of 	Problem \ref{problem: unknown length}. The TED-STG $ \Gted = (V, E, x_0) $ is a weighted directed graph constructed by:
	\begin{itemize}
		\item $ V = \cup_{t=0}^{Z} V_t $, where $ V_t = \cR_t(x_0),  \forall t \in [0, Z-1] $, $ V_{Z} = \{ \dN^0  \} $, and $ Z \deq |\cR(x_0)| $ is the reachable set size.
		\item $ E = \cup_{t=0}^{Z-2} E_t \cup E^0$, where $ E_t = \{(\delta_N^i, \delta_N^j) | \delta_N^i \in V_t,  \delta_N^j \in V_{t+1} \cap \mathcal{R}_1(\delta_N^i)    \}, \forall t \in [0, Z-2]$, and connect the terminal states in each layer to the pseudo-state $ \dN^0 $ by $ E^0 = \{ (\dNt^i, \dN^0) | \dN^i \in  \Omega \cap V_t, t \in [0, Z - 1] \} $.
	\end{itemize}
	A vertex $ \dNt^i $ above refers to the state $ \dN^i $ at time $ t $.
	The weight of the edge $ (\dNt^i, \delta_{N, t+1}^j) \in E_t, \forall t \in [0, Z-2], $ is 
	\begin{equation}
	w(\delta_{N, t}^{i}, \delta_{N, t+1}^{j}) = \min_{\dM^k \in U^{ij}} g(\delta_N^{i}, \dM^k, t),
	\end{equation}
	and the control enabling the transition from $ \delta_{N}^{i} $ to $\delta_{N}^{j}  $ at time $ t $ is also $ u_t^{ij} $ in \eqref{eq: u_t^ij}. The weight of edges in $ E^0 $ is
	\begin{equation}
		w(\dNt^i, \dN^0) = h(\dN^i, t), \ (\dNt^i, \dN^0) \in E^0.
	\end{equation}
\end{definition}

\begin{remark}
	Since the destination (terminal) states may be reached at any time in an optimal trajectory, we package all such possibilities into $ E^0 $. This greatly simplifies the problem: we only need to find an optimal path from $ \delta_{N, 0}^{i_0} $ to $ \dN^0 $. Accordingly, the time to reach the pseudo-state $ \dN^0 $ is not known in advance, and that is why it has no time subscript.
\end{remark}

\begin{figure}[tb]
	\centering
	\includegraphics[width=85mm]{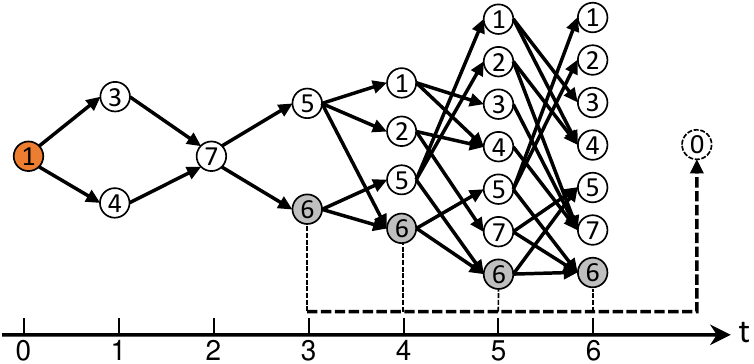} 
	\caption{ TED-STG  in Example \ref{example: TED-STG} with an initial state $x_0 = \delta_{8}^1 $ and a destination state $  \delta_8^6 $. The weight and the control associated with each edge are not shown here for clarity. Note that the pseudo-state $ \delta_8^0 $ has no fixed timestamp.}
	\label{fig: TED-STG of example}
\end{figure}

\begin{example} \label{example: TED-STG}
	We use the BCN \eqref{eq: example 1 bcn} in Example \ref{example: 1} to illustrate the TED-STG. Problem \ref{problem: unknown length} is set up by $x_0 = \{ \delta_8^1 \} $ and $ \Omega = \{ \delta_8^6 \} $. The time-variant stage cost and terminal cost are:
	\begin{equation*}
	g(x(t), u(t), t) = u(t)^{\top}Qu(t) + t, \ h(x(t), t) = x(t)^{\top}Rx(t),
	\end{equation*}
	where $ Q = \textrm{diag}(2, 3, 1, 5+t) $ and $ R = \textrm{diag}(3, 2t, 4, 0, 1, 5 + t, 6, 0) $.  The TED-STG of this task is shown in Fig. \ref{fig: TED-STG of example}. Note that the size of $ \cR(x_0) $ is $ Z = 7 $ here. We consider only one destination state  to facilitate illustration.
	\qed
\end{example}

Just like the TET-STG, the edge weights in a TED-STG do not change with time,  though the cost functions $ g $ and $ h $  are themselves time-dependent. The cost of a state trajectory is related with the weight of a path in the TED-STG as follows.

\begin{lemma} \label{lemma: path p2c2}
	Consider Case 2 of Problem \ref{problem: unknown length} under Assumption \ref{ass: nonnegative}. Given any control sequence $ u = (\dM^{k_0}, \dM^{k_1}, \cdots, \dM^{k_{K-1}} ) $ of length $ K < |\cR(x_0)|$ that steers the BCN along a state trajectory $ s = (\dN^{i_0}, \dN^{i_1}, \cdots, \dN^{i_K} ) $ with $x_0 = \delta_{N}^{i_0}  $ and $ \dN^{i_K} \in \Omega.$ There holds $ J(u)  \ge w(p) $, where $ p $ is the corresponding path from $ \delta_{N, 0}^{i_0} $ to $ \dN^0 $ in the associated TED-STG,
	\begin{equation} \label{eq: a path p2c2}
	p = \big< \delta_{N, 0}^{i_0}, \delta_{N, 1}^{i_1}, \cdots, \delta_{N, K}^{i_K}, \delta_{N}^{0}  \big>,
	\end{equation}
	and $ J(u)  = w(p) $ holds if $ \dM^{k_t} = u_t^{i_t i_{t+1}}, \forall t \in [0, K - 1]$.
\end{lemma}

Lemma \ref{lemma: path p2c2} can be proved easily in almost the same way as Lemma \ref{lemma: trajectory & path}, which is omitted here.
Now we are ready to solve Case 2 of Problem \ref{problem: unknown length} by  converting it to  an SP problem in the TED-STG through the following theorem.

\begin{theorem}\label{thm: p2 case 2}
	Consider Case 2 of Problem \ref{problem: unknown length} under Assumption \ref{ass: nonnegative} with time-dependent $ g $ and $ h $, and suppose it is feasible. Given the initial state $x_0=\delta_N^{i_0}$,  there exists an SP $ p^* = \big< \delta_{N,0}^{i_0}, \delta_{N,1}^{i_1^*}, \cdots, \delta_{N,K-1}^{i_{K-1}^*}, \delta_{N, K}^{i_K^*}, \delta_N^0 \big>, K < |\cR(x_0)|, $ in the TED-STG $ \Gted = (V, E, x_0) $  such that the minimum value of the cost function \eqref{eq: general unknown length problem} is $ J^* = w(p^*) $. The corresponding optimal control sequence is $ u^* = \{ u_t^{i_0i_1^*},  u_t^{i_1^*i_2^*}, \cdots, u_t^{i_{K-1}^* i_K^*} \} $.
\end{theorem}
\begin{IEEEproof}
	Since the problem is feasible, Proposition \ref{prop: p2} implies there exists an optimal control sequence shorter than $ |\cR(x_0)| $. Thus, we can search $ \bbU' = \{  u | |u| < |\cR(x_0)|, u\in  \bbU   \}  $ for an optimal one, where $ \bbU $ is the feasible set of Problem \ref{problem: unknown length}. Now recall Lemma \ref{lemma: path p2c2}: given any $ u \in \bbU' $, there exists a path $ p $ from $ \delta_{N, 0}^{i_0} $ to $ \dN^0 $ in the TED-STG such that $ J(u) \ge w(p) \ge w(p^*)$. Besides, with the weight given by Definition \ref{def: TED-STG}, it is obvious that $ w(p^*) = J(u^*) $. Thus, the proof is complete.
\end{IEEEproof}

Theorem \ref{thm: p2 case 2} transforms Problem \ref{problem: unknown length} with time-dependent cost functions into a standard single-pair (i.e., $ \delta_{N, 0}^{i_0} $ to $ \delta_{N}^0 $) SP problem in the TED-STG. As discussed in Remark \ref{rmk: h nonngeative}, we assume that $ h \ge 0 $ without loss of generality. Since all edges of this graph have nonnegative weights, we can apply Dijkstra's algorithm again to identify the SP, which is the same as Algorithm \ref{alg: dijkstra} except that the TED-STG is used instead of the STG$ ^+ $. We detail this algorithm in the online supplementary material on ArXiv\footnote{Refer to \url{https://arxiv.org/abs/1908.02019} for the supplementary material.} (Algorithm \ref{S-alg: dijkstra case2}) to conserve space here.

\textit{Time complexity analysis of Algorithm \ref{S-alg: dijkstra case2}.} Definition \ref{def: TED-STG} tells that there are $ Z \deq |\cR(x_0)| $ layers in the TED-STG. Each layer has at most $ Z $ vertices, and at most $ MZ $ edges exist between two consecutive layers. Besides, there are typically only few destination states, i.e., the number of incoming edges of the pseudo-state $ \delta_N^0 $ is $ O(Z) $. We thus have $ |V| =  O(Z^2)$ and $ |E| = O(MZ^2) $. Like a TET-STG, the TED-STG can be built quickly in linear time, i.e., $ O(|V| + |E|) = O(MZ^2) $. Thus, the time complexity of Algorithm \ref{S-alg: dijkstra case2} is dominated by the Dijkstra's SP part. To conclude, Algorithm \ref{S-alg: dijkstra case2} runs in time $ O(MZ^2 + Z^2\log Z^2) = O(Z^2(M + 2\log Z)) $, which is much faster than the naive Procedure \ref{proc: case 2 every n}. Since we have $ Z \le N $, the time complexity is equivalent to $ O(N^2(M + 2\log N)) $.

\begin{example}
	Recall Example \ref{example: TED-STG} and its TED-STG in Fig. \ref{fig: TED-STG of example}. Algorithm \ref{S-alg: dijkstra case2} yields the following results: $ J^* = 14 $ and $ u^* = (\delta_4^3, \delta_4^3, \delta_4^3, \delta_4^1) $. The corresponding SP in the TED-STG is $ p^* = (\delta_{8, 0}^1, \delta_{8, 1}^4, \delta_{8, 2}^7, \delta_{8, 3}^5, \delta_{8, 4}^6, \delta_{8}^0) $.
	\qed
\end{example}

\begin{table*}[h]
	\centering
	\caption{Time complexity comparison between existing work and the proposed approach on finite-horizon optimal control}
	\label{tab: time complexity}
	\begin{threeparttable}
		
		\begin{tabular}{|c|l|l|c|}
			\hline
			\multirow{2}{*}{Problem type}                          & \multicolumn{1}{c|}{\multirow{2}{*}{Task characteristics}}                                                           & \multicolumn{2}{c|}{Time complexity}                                                                                                                                                                                                                                               \\ \cline{3-4} 
			& \multicolumn{1}{c|}{}                                                                                                & \multicolumn{1}{c|}{Existing work}                                                                                                                                & Proposed approach                                                                                              \\ \hline
			\multirow{5}{*}{Problem \ref{problem: fixed length}}   & \begin{tabular}[c]{@{}l@{}}Mayer-type optimal control \\ (only terminal cost)\end{tabular}                           & $O(TMN^2)$ \cite{laschov2010maximum,laschov2013pontryagin}                                                                                                        & \multirow{5}{*}{\begin{tabular}[c]{@{}c@{}}$O(TMN)$\\ (Algorithm \ref{alg: fixed horizon})\end{tabular}}       \\ \cline{2-3}
			& \begin{tabular}[c]{@{}l@{}}Minimum-energy control \\ (only time-invariant stage cost)\end{tabular}          & $O(TN^3)$ \cite[Algorithm 3.2]{li2013minimum}                                                                                                                     &                                                                                                                \\ \cline{2-3}
			& \begin{tabular}[c]{@{}l@{}}Time-discounted stage cost \\ (no terminal cost)\end{tabular}                             & \begin{tabular}[c]{@{}l@{}}$O(MN + N^3\log_2 T)$ \cite[Theorem 3]{zhu2018optimal}\\ $O(M^2N^3 + TMN(N+M))$ \cite[Proposition 4.1]{cheng2015receding}\end{tabular} &                                                                                                                \\ \cline{2-3}
			& \begin{tabular}[c]{@{}l@{}}Both terminal  cost and stage cost\\ (general form, time-invariant)\end{tabular}                        & \begin{tabular}[c]{@{}l@{}}$O(MN + N^3\log_2 T)$ \cite[Corollary 2]{zhu2018optimal}\\ $O(TMN + TN\log(TN))$ \, \cite[Theorem 2.14]{cui2018optimal}\end{tabular}      &                                                                                                                \\ \cline{2-3}
			& \begin{tabular}[c]{@{}l@{}}Both terminal  cost and stage cost\\ (general form, possibly time-variant)\end{tabular}                 & $O(TMN)$ \cite{fornasini2013optimal}                                                                                                                              &                                                                                                                \\ \hline
			\multirow{4}{*}{Problem \ref{problem: unknown length}} & \begin{tabular}[c]{@{}l@{}}Time-optimal control \\ (only constant stage cost equal to 1)\end{tabular}                         & $O(MN^3)$\tnote{$\dagger$}  \,  \cite{laschov2013minimum, chen2016minimum}                                                                                                    & \multirow{3}{*}{\begin{tabular}[c]{@{}c@{}}$ O(MN + N\log N) $\\ (Algorithm \ref{alg: dijkstra})\end{tabular}} \\ \cline{2-3}
			& \begin{tabular}[c]{@{}l@{}}Minimum-energy control \\ (time-invariant stage cost, time not fixed)\end{tabular} & $O(N^4)$ \cite[Algorithm 3.3 ]{li2013minimum}                                                                                                                     &                                                                                                                \\ \cline{2-3}
			& \begin{tabular}[c]{@{}l@{}}Time-invariant stage cost \\ (general form, no terminal cost)\end{tabular}                                  & $ O(MN + N\log N) $ \, \cite[Theorem 2.7]{cui2018optimal}                                                                                                   &                                                                                                                \\ \cline{2-4} 
			& \begin{tabular}[c]{@{}l@{}}Time-variant stage cost\\ and (or) terminal cost\end{tabular}                               & \multicolumn{1}{c|}{---}                                                                                                                                          & \begin{tabular}[c]{@{}c@{}}$O(N^2(M + 2\log N))$\\ (Algorithm \ref{S-alg: dijkstra case2})\end{tabular}        \\ \hline
		\end{tabular}
		
		\begin{tablenotes}
			\item[$\dagger$] It is $ O(T^*MN^2) $ more precisely, where $ T^* $ is the minimum time actually required. Note that we have $ T^* = N-1$ in the worst case. 
		\end{tablenotes}
	\end{threeparttable}
\end{table*}

\section{Comparison with Related Work} \label{sec: time complexity comparison}
As we have reviewed in Section \ref{sec: intro}, unlike our algorithms which target the most general problems, most existing methods are developed for certain special cases of Problem \ref{problem: fixed length} or Problem \ref{problem: unknown length}. We therefore categorize various optimal control tasks according to their characteristics to facilitate comparison. Their time complexity is summarized in Table \ref{tab: time complexity}, where, as always, $ N=2^n $ and $ M=2^m$  for a $ n $-state, $ m $-input BCN, and $ T $ denotes the fixed horizon length in Problem \ref{problem: fixed length}.

To better interpret Table \ref{tab: time complexity}, note that we can always assume $ M \le N $ because a state can transit to at most $ N $ succeeding states regardless of the number of control inputs. In fact, we usually have $ m < n $  and thus $ M \ll N $ in practice especially for large networks. For example, it can be enough to steer the whole network by controlling only a fraction of the nodes \cite{kim2013discovery,lu2015pinning}. In short, Table \ref{tab: time complexity} shows that our graph-theoretical approach can accomplish higher time efficiency than most existing approaches, and only methods in \cite{fornasini2013optimal} and \cite{cui2018optimal} share the same time complexity as ours for particular problems.
Notably, if there are time-variant costs in Problem \ref{problem: unknown length},  only Algorithm \ref{S-alg: dijkstra case2} can handle it to the best of our knowledge. In summary, though we target FHOC problems ambitiously in their most general form,  the computational efficiency of our approach is still superior to that of most existing work. 

Note that the time complexity listed in Table \ref{tab: time complexity} refers to the \textit{worst-case} one, which indicates the longest running time of an algorithm given any possible input. By convention, the worst-case running time is used to measure the efficiency of algorithms \cite{cormen2009introduction,zhu2018optimal, wu2019optimal}. A noteworthy point  is that all the algebraic approaches in Table \ref{tab: time complexity}, i.e., all existent work except \cite{cui2018optimal}, have their average-case time complexity equal to the worst-case one, because they essentially operate on matrices of identical sizes irrespective of the initial state $ x_0 $ and the size of its reachable set $ \cR(x_0) $. By contrast, as shown in time complexity analysis of our algorithms, the actual size of the graph depends on $ \mathcal{R}(x_0) $, which is typically a small subset of the state space, while the algebraic methods always consider all the $ N $ states. Additionally, in constraint handling, our approach excludes the undesirable states and transitions completely from the graph, but most algebraic approaches simply assign them infinitely large cost values and still involve them in subsequent operations.
Consequently, our graph-theoretical approach attains potentially lower average-case time complexity than algebraic methods in practice, like \cite{fornasini2013optimal}, even though they share the same worst-case complexity. 

\section{A Biological Benchmark Example}  \label{sec: benchmark}
In this section, we focus on a larger BCN, i.e., the Ara operon gene regulatory network in \textit{E. coli} that is responsible for sugar metabolism regulation \cite{jenkins2017bistability, wu2019optimal}.  The main purpose is to compare the computational efficiency of various approaches. The Boolean model of this network is given in Table \ref{tbl: benchmark time}, where the target nodes indicate the 9 state variables, and the 4 control inputs are $ \{A_e, A_{em}, A_{ra_{-}}, G_e  \} $. Interested readers may consult \cite{jenkins2017bistability} for the biological meaning of these variables. Using the STP, we can get the ASSR of this network with a network transition matrix $ L \in \cL_{N \times MN }, N = 512, M = 16 $. 
\begin{table}
	\centering
	\caption{Boolean Model of the Ara Operon Network}
	\label{tbl: benchmark time}
	\begin{tabular}{ll} 
		\toprule
		Target node          & Boolean update rule             \\ 
		\midrule
		$A$         & $A_e \land T$                   \\
		${A_m}$     & $(A_{em}\land T)\lor A_e$       \\
		${A_{ra_+}}$ & $(A_m \lor A) \land A_{ra_{-}}$  \\
			$C$         & $\lnot G_e$                     \\
			$E$         & $M_s$                           \\
			$D $	& $ \lnot A_{ra_+} \land A_{ra_{-}}$\\
			$ M_S $ & $ A_{ra_+} \land C \land \lnot D $ \\
			$ M_T $ & $ A_{ra_+} \land C $ \\
			$ T $ & $ M_T $\\
			\bottomrule
		\end{tabular}
	\end{table}

To facilitate comparison with existing studies, we consider two tasks widely studied in the literature, i.e., the minimum-energy control and the minimum-time control. In both tasks, no constraints are enforced because only few existing methods consider state or input constraints.
Suppose the initial state is $ x_0 = \delta_{512}^{9} $ and the desired state is $ x_d = \delta_{512}^{410} $ for both tasks. All algorithms were implemented with Python 3.7. We did experiments on a desktop PC equipped with a 3.4 GHz Core i7-3770 CPU, 16 GB RAM, and Windows 10.

\begin{table*}[tb]
	\centering
	\caption{Running Time of Various Methods for Optimal Control of the Ara Operon Network}
	\label{tbl: running time}
	\begin{tabular}{c|ccccc|cccc} 
		\toprule
		Problem                               & \multicolumn{5}{c|}{Task 1}                                            & \multicolumn{4}{c}{Task 2}  \\ 
		\midrule
		Method                                & \cite[Algorithm 3.2]{li2013minimum}             & \cite{fornasini2013optimal}            & \cite{zhu2018optimal}              & \cite[Theorem 2.14]{cui2018optimal}         & Ours & \cite[Algorithm 3.3]{li2013minimum} & \cite{laschov2013minimum} & \cite[Theorem 2.7]{cui2018optimal} & Ours         \\
		\multicolumn{1}{l|}{Running time (s)} & 350.43& 0.23 & 402.79 & 0.17 & 0.10                   & 19651.33   &  636.72  &  0.01  &     0.01         \\
		\bottomrule
	\end{tabular}
\end{table*}

\subsection{Task 1: Minimum-Energy Control}
In this task, we aim to transfer the BCN from the initial state $ x_0 $ to the desired state $ x_d $ at a prespecified time point with least energy consumption \cite{li2013minimum}. This task is presented as an instance of Problem \ref{problem: fixed length} with $ T = 10 $ and $ \Omega = \{ x_d \} $. Since most methods can only deal with time-invariant costs, we use a stage cost function like that in \cite{wu2019optimal} to evaluate the energy consumed by each state transition, and set zero terminal costs:
\begin{equation*}
g(x(t), u(t), t) = g(\ltimes_{i=1}^9x_i(t), \ltimes_{j=1}^4u_j(t)) = \mathcal{A}X(t) + \mathcal{B}U(t),
\end{equation*}
where $ X(t) = [x_1(t), x_2(t), \cdots, x_9(t)]^{\top}$ and $ U(t) = [u_1(t), u_2(t), u_3(t), u_4(t)]^{\top}$. The two weight vectors are  $ \mathcal{A} = [0, 16, 40, 44, 28, 28, 28, 48, 44] $, and $ \mathcal{B} =  [0, 48, 28, 48]$.

 Running Algorithm \ref{alg: fixed horizon}, we get the following results.
\begin{itemize}
	\item The reachable set of $ x_0 $, i.e., $ \cR(x_0) $, has only 108 states, though the complete state state has 512 states in total.
	\item The minimum value of $ \JT(u) $ in \eqref{eq: general fixed length problem} is $ \JT^* = 1108 $.
	\item The optimal control sequence is $ u^* = \delta_{16}(16, 16, 16, 16, 16, 16, 8, 5, 6, 14) $,
	and the resultant state trajectory is $ 	s^* = \delta_{512}(9, 457, 463, 480, 480, 480, 480, 352, 312, 288, 410). $
\end{itemize}

\begin{figure}[tb]
	\centering
	\includegraphics[width=87mm]{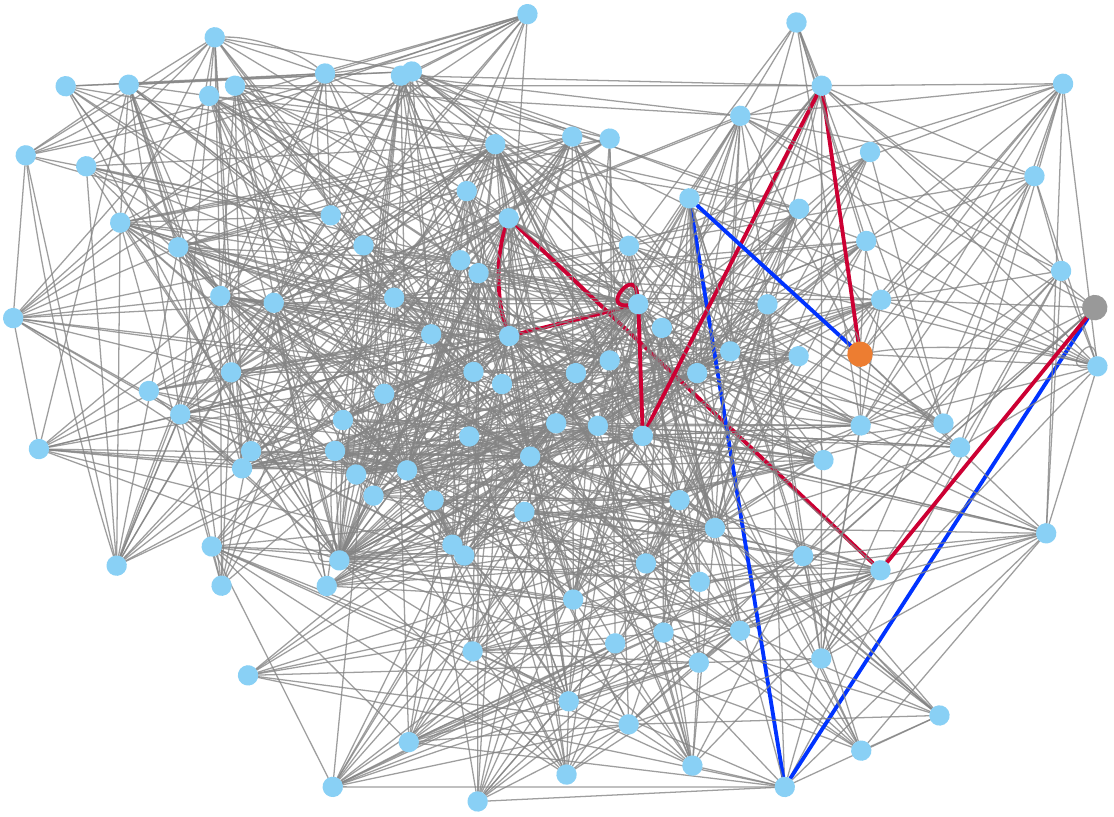} 
	\caption{ The state transition graph of the Ara operon network with the initial state $ x_0 = \delta_{512}^9 $ (in orange) and the desired state $ x_d = \delta_{512}^{410}  $ (in gray). The state trajectories of minimum-energy and minimum-time control are highlighted in red and blue respectively.  The direction of each state transition and the label of each state are not shown here for readability.}
	\label{fig: benchmark}
\end{figure}

We illustrate $ \cR(x_0) $ and highlight the above optimal state trajectory $ s^* $ in Fig. \ref{fig: benchmark}. We tested the other methods, and they all obtained the same optimal value $ J_T^* $. The running time of different methods is listed in Table \ref{tbl: running time}. As we have expected, the algorithm in \cite{fornasini2013optimal} takes more time than ours in practice, though they have identical worst-case time complexity. The main reason is that the former always evaluates the whole state space, while our algorithm only focuses on the reachable states $ \cR(x_0) $.  Another method \cite{cui2018optimal} shares the same view as ours, but it depends on the slower Dijkstra's algorithm rather than the more efficient DP-based Algorithm \ref{alg: fixed horizon}. Overall, our approach attains the shortest running time to complete Task 1.

\subsection{Task 2: Minimum-Time Control}
This task is also referred to as time-optimal control,  that is, to steer the BCN from $ x_0 $ to a desired state $ x_d $ as fast as possible \cite{laschov2013minimum}. It is easy to specialize Problem \ref{problem: unknown length} for this task: just let $ g(\cdot, \cdot, \cdot) \equiv 1 $, i.e., unit time for each state transition, $ h(\cdot, \cdot) \equiv 0 $, and $ \Omega = \{x_d\} $. In this setting, the optimal value $ J^* $ of \eqref{eq: general unknown length problem} is obviously the minimum time taken from $ x_0 $ to $ x_d $.
Since $ g $ and $ h $ are both time-independent, we can tackle this task with Algorithm \ref{alg: dijkstra}. The results are given below.
\begin{itemize}
	\item The optimal value is $ J^* = 3 $, i.e., the BCN is transferred from $ x_0 $ to $x_d  $ in at least 3 steps.
	\item The optimal control sequence is $ u^* = \delta_{16}(1, 2, 14) $, which leads to the state trajectory $ s^* = \delta_{512}(9, 41, 15, 410) $.
\end{itemize}

Like Task 1, we sketch the above minimum-time state trajectory in Fig. \ref{fig: benchmark}. 
The other methods yielded the same minimum time, and their running time is listed in Table \ref{tbl: running time}.  In this task, the method in \cite{cui2018optimal} is essentially identical to ours, both depending on Dijkstra's algorithm to find an SP, and they are the fastest ones, taking far less time than the other two. Note that we have optimized it by building the STG efficiently with Algorithm \ref{alg: STG} instead of the expensive algebraic method originally used in \cite{cui2018optimal}. The algorithm in \cite{li2013minimum} is extremely slow here mainly because it examines exhaustively all SP's of length ranging from 1 to $ N $ to find the shortest one.

\begin{remark}
	The running time of some algebraic approaches, like \cite{zhu2018optimal} and \cite{laschov2013minimum}, might be further reduced using advanced numerical routines, since they depend heavily on matrix operations, and the involved matrices are typically sparse. Nevertheless, the results in Table \ref{tbl: running time} still demonstrate the supreme efficiency of our approach with an advantage of several orders of magnitude. Besides, the method in \cite{cui2018optimal} is closest to ours, but it can only solve a small subset of problems investigated in this study. The Python implementation of the proposed approach and existing algorithms is available at GitHub \url{https://github.com/ShuhuaGao/FHOC}.
\end{remark}

\section{Conclusion} \label{sec: conclusion}
This paper focused on FHOC of BCNs from a graph-theoretical perspective. We unified various kinds of specific FHOC problems into two general constrained optimization problems, which can incorporate time-variant costs and a diverse range of constraints. Then, as a central contribution of this study, we established equivalence between general FHOC problems and the SP problem in particular graphs. Two efficient algorithms were afterwards designed to find such an SP.  As shown by both time complexity analysis and numerical experiments, our approach can handle the most general problems while maintaining a competitive advantage in computational efficiency. Finally, we note that all SP problems in Problem \ref{problem: fixed length} and \ref{problem: unknown length} can be technically solved by a single SP algorithm, like Dijkstra's algorithm, though we proposed two custom algorithms for efficiency purpose. That's why we consider our graph-theoretical approach as a unified framework, which is characterized by high computational efficiency and methodological consistency across a wide range of FHOC problems.
Due to the discrete and deterministic nature of BCNs, we believe it is a promising research direction to hybridize the newly developed ASSR with the classical graph theory for more studies on BCNs beyond FHOC. One future work is to adapt this graph-theoretical approach to infinite-horizon optimal control problems.


%

\ifCLASSOPTIONcaptionsoff
  \newpage
\fi



%
\bibliographystyle{IEEEtran}
\bibliography{boolnet}

\begin{thebibliography}{10}
\providecommand{\url}[1]{#1}
\csname url@samestyle\endcsname
\providecommand{\newblock}{\relax}
\providecommand{\bibinfo}[2]{#2}
\providecommand{\BIBentrySTDinterwordspacing}{\spaceskip=0pt\relax}
\providecommand{\BIBentryALTinterwordstretchfactor}{4}
\providecommand{\BIBentryALTinterwordspacing}{\spaceskip=\fontdimen2\font plus
\BIBentryALTinterwordstretchfactor\fontdimen3\font minus
  \fontdimen4\font\relax}
\providecommand{\BIBforeignlanguage}[2]{{%
\expandafter\ifx\csname l@#1\endcsname\relax
\typeout{** WARNING: IEEEtran.bst: No hyphenation pattern has been}%
\typeout{** loaded for the language `#1'. Using the pattern for}%
\typeout{** the default language instead.}%
\else
\language=\csname l@#1\endcsname
\fi
#2}}
\providecommand{\BIBdecl}{\relax}
\BIBdecl

\bibitem{kauffman1969metabolic}
S.~A. Kauffman, ``Metabolic stability and epigenesis in randomly constructed
  genetic nets,'' \emph{Journal of theoretical biology}, vol.~22, no.~3, pp.
  437--467, 1969.

\bibitem{saadatpour2013boolean}
A.~Saadatpour and R.~Albert, ``Boolean modeling of biological regulatory
  networks: a methodology tutorial,'' \emph{Methods}, vol.~62, no.~1, pp.
  3--12, 2013.

\bibitem{datta2003external}
A.~Datta, A.~Choudhary, M.~L. Bittner, and E.~R. Dougherty, ``External control
  in markovian genetic regulatory networks,'' \emph{Machine learning}, vol.~52,
  no. 1-2, pp. 169--191, 2003.

\bibitem{caetano2015boolean}
M.~A.~L. Caetano and T.~Yoneyama, ``Boolean network representation of contagion
  dynamics during a financial crisis,'' \emph{Physica A: Statistical Mechanics
  and its Applications}, vol. 417, pp. 1--6, 2015.

\bibitem{cheng2009controllability}
D.~Cheng and H.~Qi, ``Controllability and observability of boolean control
  networks,'' \emph{Automatica}, vol.~45, no.~7, pp. 1659--1667, 2009.

\bibitem{cheng2010linear}
D.~Cheng and H.~Qi, ``A linear representation of dynamics of boolean
  networks,'' \emph{IEEE Transactions on Automatic Control}, vol.~55, no.~10,
  pp. 2251--2258, 2010.

\bibitem{zhao2010input}
Y.~Zhao, H.~Qi, and D.~Cheng, ``Input-state incidence matrix of boolean control
  networks and its applications,'' \emph{Systems \& Control Letters}, vol.~59,
  no.~12, pp. 767--774, 2010.

\bibitem{laschov2013observability}
D.~Laschov, M.~Margaliot, and G.~Even, ``Observability of boolean networks: A
  graph-theoretic approach,'' \emph{Automatica}, vol.~49, no.~8, pp.
  2351--2362, 2013.

\bibitem{liang2017algorithms}
J.~Liang, H.~Chen, and Y.~Liu, ``On algorithms for state feedback stabilization
  of boolean control networks,'' \emph{Automatica}, vol.~84, pp. 10--16, 2017.

\bibitem{li2019robustness}
B.~Li, Y.~Liu, J.~Lou, J.~Lu, and J.~Cao, ``The robustness of outputs with
  respect to disturbances for boolean control networks,'' \emph{IEEE
  transactions on neural networks and learning systems}, 2019.

\bibitem{zhao2015control}
Y.~Zhao, B.~K. Ghosh, and D.~Cheng, ``Control of large-scale boolean networks
  via network aggregation,'' \emph{IEEE transactions on neural networks and
  learning systems}, vol.~27, no.~7, pp. 1527--1536, 2015.

\bibitem{laschov2010maximum}
D.~Laschov and M.~Margaliot, ``A maximum principle for single-input boolean
  control networks,'' \emph{IEEE Transactions on Automatic Control}, vol.~56,
  no.~4, pp. 913--917, 2010.

\bibitem{laschov2013pontryagin}
D.~Laschov and M.~Margaliot, ``A pontryagin maximum principle for multi-input
  boolean control networks,'' \emph{Recent advances in dynamics and control of
  neural networks}, 2013.

\bibitem{li2013minimum}
F.~Li and X.~Lu, ``Minimum energy control and optimal-satisfactory control of
  boolean control network,'' \emph{Physics Letters A}, vol. 377, no.~43, pp.
  3112--3118, 2013.

\bibitem{laschov2013minimum}
D.~Laschov and M.~Margaliot, ``Minimum-time control of boolean networks,''
  \emph{SIAM Journal on Control and Optimization}, vol.~51, no.~4, pp.
  2869--2892, 2013.

\bibitem{fornasini2013optimal}
E.~Fornasini and M.~E. Valcher, ``Optimal control of boolean control
  networks,'' \emph{IEEE Transactions on Automatic Control}, vol.~59, no.~5,
  pp. 1258--1270, 2013.

\bibitem{zhu2018optimal}
Q.~Zhu, Y.~Liu, J.~Lu, and J.~Cao, ``On the optimal control of boolean control
  networks,'' \emph{SIAM Journal on Control and Optimization}, vol.~56, no.~2,
  pp. 1321--1341, 2018.

\bibitem{cheng2015receding}
D.~Cheng, Y.~Zhao, and T.~Xu, ``Receding horizon based feedback optimization
  for mix-valued logical networks,'' \emph{IEEE Transactions on Automatic
  Control}, vol.~60, no.~12, pp. 3362--3366, 2015.

\bibitem{zhao2010optimal}
Y.~Zhao, Z.~Li, and D.~Cheng, ``Optimal control of logical control networks,''
  \emph{IEEE Transactions on Automatic Control}, vol.~56, no.~8, pp.
  1766--1776, 2011.

\bibitem{zhao2011floyd}
Y.~Zhao, ``A floyd-like algorithm for optimization of mix-valued logical
  control networks,'' in \emph{Proceedings of the 30th Chinese Control
  Conference}.\hskip 1em plus 0.5em minus 0.4em\relax IEEE, 2011, pp.
  1972--1977.

\bibitem{cheng2014optimal}
D.~Cheng, Y.~Zhao, and J.-B. Liu, ``Optimal control of finite-valued
  networks,'' \emph{Asian Journal of Control}, vol.~16, no.~4, pp. 1179--1190,
  2014.

\bibitem{wu2019optimal}
Y.~Wu, X.-M. Sun, X.~Zhao, and T.~Shen, ``Optimal control of boolean control
  networks with average cost: A policy iteration approach,'' \emph{Automatica},
  vol. 100, pp. 378--387, 2019.

\bibitem{cui2018optimal}
X.~Cui, J.-E. Feng, and S.~Wang, ``Optimal control problem of boolean control
  networks: A graph-theoretical approach,'' in \emph{2018 Chinese Control And
  Decision Conference (CCDC)}.\hskip 1em plus 0.5em minus 0.4em\relax IEEE,
  2018, pp. 4511--4516.

\bibitem{zhang2017finite}
Z.~Zhang, T.~Leifeld, and P.~Zhang, ``Finite horizon tracking control of
  boolean control networks,'' \emph{IEEE Transactions on Automatic Control},
  vol.~63, no.~6, pp. 1798--1805, 2017.

\bibitem{faryabi2008optimal}
B.~Faryabi, G.~Vahedi, J.-F. Chamberland, A.~Datta, and E.~R. Dougherty,
  ``Optimal constrained stationary intervention in gene regulatory networks,''
  \emph{EURASIP Journal on Bioinformatics and Systems Biology}, vol. 2008,
  no.~1, p. 620767, 2008.

\bibitem{akutsu2007control}
T.~Akutsu, M.~Hayashida, W.-K. Ching, and M.~K. Ng, ``Control of boolean
  networks: Hardness results and algorithms for tree structured networks,''
  \emph{Journal of theoretical biology}, vol. 244, no.~4, pp. 670--679, 2007.

\bibitem{liang2017improved}
J.~Liang, H.~Chen, and J.~Lam, ``An improved criterion for controllability of
  boolean control networks,'' \emph{IEEE Transactions on Automatic Control},
  vol.~62, no.~11, pp. 6012--6018, 2017.

\bibitem{zhu2018further}
Q.~Zhu, Y.~Liu, J.~Lu, and J.~Cao, ``Further results on the controllability of
  boolean control networks,'' \emph{IEEE Transactions on Automatic Control},
  vol.~64, no.~1, pp. 440--442, 2018.

\bibitem{pal2006optimal}
R.~Pal, A.~Datta, and E.~R. Dougherty, ``Optimal infinite-horizon control for
  probabilistic boolean networks,'' \emph{IEEE Transactions on Signal
  Processing}, vol.~54, no.~6, pp. 2375--2387, 2006.

\bibitem{kim2013discovery}
J.~Kim, S.-M. Park, and K.-H. Cho, ``Discovery of a kernel for controlling
  biomolecular regulatory networks,'' \emph{Scientific reports}, vol.~3, p.
  2223, 2013.

\bibitem{o1981discrete}
J.~O'Reilly, ``The discrete linear time invariant time-optimal control
  problem—an overview,'' \emph{Automatica}, vol.~17, no.~2, pp. 363--370,
  1981.

\bibitem{chen2016minimum}
H.~Chen, B.~Wu, and J.~Lu, ``A minimum-time control for boolean control
  networks with impulsive disturbances,'' \emph{Applied Mathematics and
  Computation}, vol. 273, pp. 477--483, 2016.

\bibitem{cormen2009introduction}
T.~H. Cormen, C.~E. Leiserson, R.~L. Rivest, and C.~Stein, \emph{Introduction
  to Algorithms}, 3rd~ed.\hskip 1em plus 0.5em minus 0.4em\relax The MIT Press,
  2009.

\bibitem{pallottino1998shortest}
S.~Pallottino and M.~G. Scutella, ``Shortest path algorithms in transportation
  models: classical and innovative aspects,'' in \emph{Equilibrium and advanced
  transportation modelling}.\hskip 1em plus 0.5em minus 0.4em\relax Springer,
  1998, pp. 245--281.

\bibitem{lu2015pinning}
J.~Lu, J.~Zhong, C.~Huang, and J.~Cao, ``On pinning controllability of boolean
  control networks,'' \emph{IEEE Transactions on Automatic Control}, vol.~61,
  no.~6, pp. 1658--1663, 2015.

\bibitem{jenkins2017bistability}
A.~Jenkins and M.~Macauley, ``Bistability and asynchrony in a boolean model of
  the l-arabinose operon in escherichia coli,'' \emph{Bulletin of mathematical
  biology}, vol.~79, no.~8, pp. 1778--1795, 2017.

\end{thebibliography}

\includepdf[pages=-, openright=true]{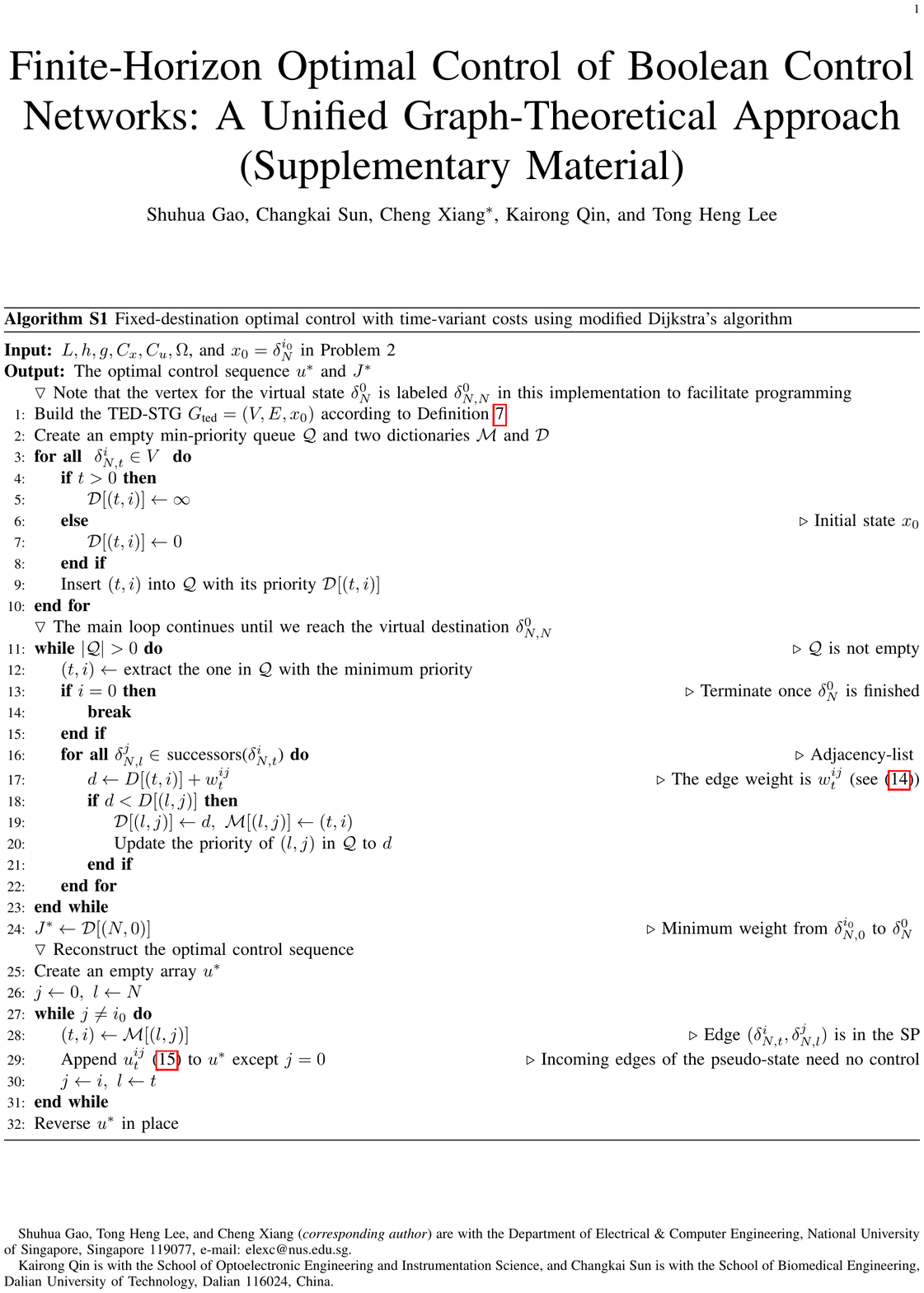}
\end{document}